\documentclass[3p,review,llncs,10pt]{elsarticle}
\usepackage{amssymb}
\usepackage{amsmath}
\usepackage{theorem}
\usepackage{graphicx}
\usepackage{subcaption}
\usepackage{multirow}
\usepackage[linesnumbered,ruled,resetcount, noend]{algorithm2e}
\usepackage[]{natbib}
\usepackage[hidelinks]{hyperref}
\usepackage{geometry}
\usepackage{rotating}
\usepackage{float}
\usepackage{booktabs,caption}
\usepackage[referable, flushleft]{threeparttablex}
\usepackage{longtable}
\usepackage{tikz}
\usepackage{xxcolor}
\usepackage{pgfplots}
\usepackage{environ}
\usepackage{subcaption}
\usepackage{setspace}
\usepackage{enumitem}
\usepackage{siunitx}

\usepackage[abs]{overpic}

\newtheorem{theo}{Theorem}[section]
{\theorembodyfont{\rmfamily}

}

{\theorembodyfont{\rmfamily}
\newtheorem{example}[theo]{Example}
}
{\theorembodyfont{\rmfamily}

}
{\theorembodyfont{\rmfamily}

}

\newcommand{\ubar}[1]{\text{\b{$#1$}}}
\newcommand{\var}[2]{{#1}_{#2}^{id}}
\newcommand{\prm}[2]{{#1}_{#2}^{i}}
\newcommand{\tw}[2]{{#1}_{#2}^{d}}

\begin{document}

\title{A biobjective Home Care Scheduling Problem with dynamic breaks}

\author[modes]{Isabel M\'endez-Fern\'andez}
\ead{isabel.mendez.fernandez@udc.es}

\author[modes]{Silvia Lorenzo-Freire}
\ead{silvia.lorenzo@udc.es}

\author[modestya,citma]{Ángel Manuel González-Rueda\corref{mycorrespondingauthor}}
\ead{angelmanuel.gonzalez.rueda@usc.es}

\address[modes]{MODES Research Group, Department of Mathematics and Centre for Information and Communications Technology Research (CITIC)\\ Faculty of Computer Science, University of A Coruña, Campus de Elviña, A Coru\~{n}a, Spain}

\address[modestya]{MODESTYA Research Group, Department of Statistics, Mathematical Analysis and Optimization and. University of Santiago de Compostela. Santiago de Compostela, Spain}

\address[citma]{CITMAga (Galician Center for Mathematical Research and Technology), 15782 Santiago de Compostela, Spain}

\cortext[mycorrespondingauthor]{Corresponding author. Postal code: 15782, Santiago de Compostela, Spain.}

\begin{highlights}
\item We present a real biobjective Home Care Scheduling Problem.
\item We propose a biobjective Mixed Integer Linear Programming formulation.
\item We design a custom metaheuristic based on the Multi-Directional Local Search (MDLS) technique.
\item We compare the MDSL algorithm with other two well known multiobjective methods: AUGMECON2 and NSGA-II.
\item An extensive computational analysis is conducted including a real case study. 

\end{highlights}

\begin{abstract}
    This paper presents a multiobjective Home Care Scheduling Problem (from now on multiobjective HCSP) related to a home care company for elderly and dependent people located in the North of Spain. In particular, a biobjective problem is considered, with the following two conflicting objectives: the welfare of users and the cost of schedules. 
    

    To tackle the problem, a custom metaheuristic algorithm based on the Multi-Directional Local Search (MDLS) was designed, obtaining good approximations of the Pareto frontier in efficient computational times. This biobjective algorithm can be divided into three steps: initializing the set of non dominated solutions, generating solutions composed by different routes and obtaining non dominated solutions.  

    The performance of the biobjective algorithm was analyzed by implementing two other well known methods in the literature: the exact method AUGMECON2, which is just an improved version of the Epsilon Constraint approach, and an NSGA-II-based algorithm. 
    

    Finally, an extensive computational study was developed to compare the three methods over a set of instances from the literature, where the biobjective algorithm exhibited a superior behaviour. Furthermore, the algorithm was also applied to real instances providing solutions to the company with a good trade-off between the two objectives.

\end{abstract}


\maketitle

\section{Introduction}

Home care is a service that allows elderly and/or dependent people to continue living in their homes despite being in a situation of dependency or in need of assistance to  undertake day-to-day tasks. This type of service improves the quality of life of users and their families, providing domestic and care services to ensure that users are able to enjoy the best possible living conditions.

A multiobjective Home Care Scheduling Problem (multiobjective HCSP) involves designing routes and schedules for caregivers, indicating the services to carry out, in which order and at what time while optimizing more than one factor. In recent years, the demand of home care services has grown, which has resulted in an increase of the literature on this type of problems.  Table~\ref{table:lit} shows the most common characteristics of the multiobjective HCSP studied in the literature. 

\renewcommand{\arraystretch}{1.3} 
\begin{table}[H]
	\resizebox{\textwidth}{!}{
		\begin{tabular}{l  c c c cc c c }
			\hline
			Reference                        & Service hard & Service soft & Caregiver   & Fixed services &   Maximum    & Compatibility & Continuity \\
			                                 & time windows   & time windows   & time windows & duration       & working time &               & of care    \\ \hline
			\citet{AitHaddadene2019}         & X             & -             & -            & X              &      X       & X             & -          \\
            \citet{Belhor2023}               & X             & X             & -            & -              &      X       & -             & -          \\ 
            \citet{BRAEKERS2016428}          & X             & X             & X            & X              &      X       & X             & -          \\
			\citet{DECERLE2019712}           & -             & X             & X            & X              &      -       & X             & -          \\
			\citet{DECERLE2018346}           & -             & X             & X            & X              &      -       & X             & -          \\
			\citet{FathollahiFard2020}       & X             & -             & -            & X              &      -       & -             & -          \\
			\citet{FATHOLLAHIFARD2019118200} & X             & -             & -            & X              &      -       & -             & -          \\
			\citet{FATHOLLAHIFARD2018423}    & X             & -             & -            & X              &      -       & -             & -          \\
            \citet{HabibnejadLedari2019}     & -             & -             & -            & -              &      X       & X             & -          \\
            \citet{Haddadene2016}            & X             & -             & -            & X              &      X       & X             & -          \\
            \citet{Khodabandeh2021}          & X             & -             & -            & X              &      -       & X             & -          \\
			\citet{Liu2018}                  & -             & -             & -            & -              &      X       & X             & X          \\
			\citet{Ma2022}                   & -             & X             & -            & X              &      X       & X             & -          \\
            \citet{Vieira2022}               & X             & -             & X            & X              &      X       & X             & X          \\
            \citet{Xiang2023}                & X             & -             & -             & X              &      X       & X             & X           \\
            \citet{Yang2021}                 & -             & X             & -            & -              &      X       & -             & X          \\
			Our approach                     & X             & X             & X            & X              &      X       & X             & X          \\ \hline
		\end{tabular}
	} \caption{Brief summary of the multiobjective HCSP characteristics in the literature.} \label{table:lit}
\end{table}

In this class of problems, two types of services can be considered: social care services (which include tasks like cleaning, cooking, etc) and health care services, more related to medical tasks. In certain problems, both types of tasks are combined. This is the case of \citet{Vieira2022}, where there are also time dependencies between services. 

In most cases, there is a period of time that establishes when each service can be performed, called service hard time window. In case it is desirable that the service is carried out over a period of time, the service soft time window is indicated. Both characteristics are available in \citet{BRAEKERS2016428}, providing two soft time windows per service: a tight one and a loose one. Since  soft time windows are not mandatory requirements, it is usual to penalize the deviation from them in the objective function. With the aim of avoiding an excessive deviation, in \citet{DECERLE2019712} and \citet{DECERLE2018346} the rate of the penalization increases in case the deviation exceeds a certain threshold. It is also very common in HCSP to consider fixed services duration, although in some HCSP the duration of the services may depend on other factors, such as the skills or the efficiency of the caregivers that will perform them (\citet{Belhor2023} , \citet{HabibnejadLedari2019}, \citet{Liu2018} and \citet{Yang2021}).

As in the case of services, caregivers may indicate a period of availability, known as caregiver time windows. In most of the problems there is a maximum working time (per day, week or month or even a combination of them) that cannot be surpassed. In the case of  \citet{FathollahiFard2020}, \citet{FATHOLLAHIFARD2019118200} and \citet{FATHOLLAHIFARD2018423}, the maximum working time is replaced with a limit for the travel distance. In turn, \citet{Belhor2023} consider a maximum daily working time and establish a maximum size for the set of patients each caregiver will attend.

To assign caregivers to services, it is important to take into account the compatibility of the caregiver with the service to carry out. In many cases, it is considered that the caregiver and the service are compatible if the caregiver possesses the skills required to successfully carry out the service. \citet{AitHaddadene2019}, \citet{Haddadene2016}, \citet{Ma2022} and \citet{Xiang2023} study problems where any caregiver can perform a service with lower requirements but services cannot be assigned to caregivers that do not reach the required skill. According to \citet{Khodabandeh2021}, the over-qualification of caregivers can lead to economic losses as caregivers are paid according to their skills. In other cases, the compatibility between services and caregivers is determined by the preferences of the users who demand the services. In \citet{BRAEKERS2016428} users may establish their preferences regarding the caregivers, which are indicated by a penalty (0, 1 or 2).

The continuity of care is a concept that also appears in some contributions related to HCSP, although the interpretation of the term also depends on the context of the problem to be tackled. In \citet{Liu2018} the continuity of care is interpreted by setting a maximum number of medical teams that can be assigned to a patient. For \citet{Vieira2022}, the continuity of care depends on the type of service: in the case of health services only one caregiver per user is allowed, whereas for social services it is desirable not to exceed two caregivers per user and week. 

The problem under study is based on a real company located in the Northwest of Spain, which provides home care services (mainly social care services) to elderly and dependent people. As it can be seen in Table~\ref{table:lit}, services are characterized by a fixed duration, a soft time window and a hard time window. Caregivers also have a time window that establishes when they are willing to work. Moreover, there is a agreed working time per week and a daily maximum working time for caregivers. The weekly time can be surpassed as long as the daily one is not exceeded and caregivers are paid for their overtime. In accordance with the company's guidelines, on the basis of its previous experience and needs, a combination of caregiver-user compatibility and continuity of care is adopted, measuring its importance in terms of six different levels, which reflect the degree of affinity between caregivers and users. 

Additionally, it is necessary to introduce a new feature into the problem, since the longest break that a caregiver takes during the working day will be deducted from her worked time, in case it reaches a certain duration. In contrast to other works (see, for instance, \citet{Vieira2022}, where part time caregivers have a 20 minutes break, if they work at least 2 consecutive hours, and full time caregivers have a meal break that typically takes 30 minutes), the duration of the break is not previously fixed and is determined by the final daily schedule of the caregiver. This dynamic break may not even exist if there is no rest between consecutive services for more than a specified duration. Moreover, the dynamic break will have a significant impact on the cost of caregivers, since it will be deducted from their daily working times. 

\renewcommand{\arraystretch}{1.3} 
\begin{table}[H]
	\resizebox{\textwidth}{!}{
		\begin{threeparttable}
			\begin{tabular}{l ccccccccc}
				\hline
                Reference & \multicolumn{5}{c}{Issue} & No. issues & No. objectives & \multicolumn{2}{c}{Solution method} \\
                \cline{2-6} \cline{9-10}
				& Overtime & Idle time & Travel time & Service soft &    User  &  &  & Exact & (Meta)heuristic \\
				&  & & & time windows  & preferences   & & & &\\ \hline
                \citet{AitHaddadene2019}         &       -       &      -      &   X    &    -    &    X    &  2 & 2 & - & NSGA-II\\ 
                \citet{Belhor2023}               &       -       &      -      &   X    &    X    &    -    &  3 & 2 & - & NSGA-II/SPEA2\\ 
                \citet{BRAEKERS2016428}          &       X       &      -      &   X    &    X    &    X    &  4 & 2 & $\epsilon$-constraint & MDLS\\
				\citet{DECERLE2019712}           &       -       &      X      &   X    &    X    &    -    &  5 & 3 & - & MDLS/NSGA-II\\
				\citet{DECERLE2018346}           &       -       &      -      &   X    &    X    &    -    &  4 & 3 & - & MDLS\\
				\citet{FathollahiFard2020}       &       -       &      -      &   X    &    -    &    X    &  6 & 2 & - & MOSA\\
				\citet{FATHOLLAHIFARD2019118200} &       -       &      -      &   X    &    -    &    -    &  8 & 2 & - & MOSA\\
				\citet{FATHOLLAHIFARD2018423}    &       -       &      -      &   X    &    -    &    -    &  3 & 2 & - & MOSA\\
				\citet{HabibnejadLedari2019}     &       -       &      -      &   -    &    -    &    X    &  4 & 3 & - & NSGA-II\\
                \citet{Haddadene2016}            &       -       &      -      &   X    &    -    &    X    &  2 & 2 & - & NSGA-II\\
                \citet{Khodabandeh2021}          &       -       &      -      &   X    &    -    &    -    &  2 & 2 & $\epsilon$-constraint & -\\
				\citet{Liu2018}                  &       X       &      -      &   X    &    -    &    -    &  4 & 2 & $\epsilon$-constraint & MOSA/NSGA-II\\
				\citet{Ma2022}                   &       -       &      -      &   -    &    X    &    -    &  2 & 2 & - & MOBSO/NSGA-II\\
                \citet{Vieira2022}               &       -       &      X      &   X    &    -    &    -    &  4 & 3 & - & Heuristic\\
                \citet{Xiang2023}                &       X       &      X      &   X    &    -    &    X    &  7 & 2  &  $\epsilon$-constraint & NSGA-II\\
                \citet{Yang2021}                 &       -       &      -      &   X    &    X    &    -    &  5 & 3 & - & MOABC\\
				Our approach                     &       X       &      X      &   X    &    X    &    X    &  5 & 2 & $\epsilon$-constraint & MDLS/NSGA-II\\ 
    \hline
			\end{tabular}
   \begin{tablenotes}
				\item MDLS = Multi-Directional Local Search, MOABC= Multi-Objective Artificial Bee Colony, MOBSO= Multi-Objective Brain Storm Optimization, MOSA= Multi-Objective Simulated Annealing, NSGA=  Non-dominated Sorting Genetic Algorithm, SPEA= Strength Pareto Evolutionary Algorithm.
			\end{tablenotes}
		\end{threeparttable}
	} 
	\caption{Brief summary of issues, objectives and solution methods in multiobjective HCSP.} \label{table:lit_objs}
\end{table}

Multiple issues get involved in the objective functions of multiobjective HCSP. It can be seen in Table~\ref{table:lit_objs}, which also includes some common issues to evaluate the quality of solutions on this type of problems.  Usually, the quality of solutions is influenced by the total travel time between consecutive services of the schedules of caregivers, since the working time of caregivers depends on it. Another relevant factor that may impact on the working time of caregivers is the idle or free time of caregivers between services. Idle times can be produced when caregivers wait to match the time windows of services or with the aim of synchronizing visits of several caregivers if they share services. When caregivers have an agreed working time that can be surpassed, the cost of schedules can be reduced by optimizing the resulting overtime. On the other hand, the satisfaction of users can be increased by taking into account the preferences of users in terms of the caregivers who could care for them or respecting the soft time windows of services as much as possible. 

Although the different multiobjective HCSP often involve common issues, their treatment in the objective functions may vary. In some cases, the authors consider the issues under study individually, but it is also possible to group related issues into a weighted sum to reduce the number of objectives to optimize (in fact, the number of objectives in multiobjective HCSP is usually no more than 3). Thus, in terms of the cost of schedules, \citet{BRAEKERS2016428} combine overtime and travel time into a single objective, whereas \citet{Ma2022} interpret it in terms of travel and service costs. In \citet{Vieira2022} the non effective working time is composed of travel and idle times. For \citet{Yang2021} the cost of schedules depends on travel, service and soft time window penalization costs. Regarding  satisfaction, \citet{BRAEKERS2016428} measure satisfaction of users in terms of respecting their preferences and soft time windows of services. Similarly, \citet{DECERLE2019712} and \citet{DECERLE2018346} consider the quality of services as a combination of soft time window and synchronization penalization. In addition to these more common objectives, \citet{FATHOLLAHIFARD2019118200} optimize the environmental impact of  schedules, which is the result of environmental pollution  and $\mathrm{CO_2}$ emissions.

As presented in Table~\ref{table:lit_objs}, the issues considered in the problem under study are: overtime, travel time, idle time, service soft time windows and preferences of users. These factors are combined into two objectives:
\begin{itemize}
	\item Minimize the cost of schedules, which is just a combination of the overtime and working time of caregivers. Since the duration of services is fixed, total working time can be minimized by reducing travel or idle times. Notice that, as explained before, the largest daily break is not considered as idle time if it reaches certain duration.
 \item Maximize the welfare of users, interpreted in terms of the preferences of users and the soft time window penalization for services. The soft time window penalization takes into account the performance of the services both before and after their time window. The preferences of users are indicated by the affinity level between services and caregivers.
\end{itemize}

The problem of obtaining the schedule for this company was also studied in \citet{MendezFernandez2020} and \citet{MendezFernandez2023}, but there are significant differences with the problem presented in this framework. In \citet{MendezFernandez2020} the goal was to solve the initial problem presented by the company, which consisted in updating the prearranged schedules of the caregivers in order to solve a set of incidents that regularly arose in the company. To solve the problem, a simulated annealing algorithm was developed to slightly modify the prearranged schedules of the services. 
\citet{MendezFernandez2023} tackled the problem  in a more general manner, considering that the goal was to obtain the best possible schedules without being restricted to modifying previous schedules. Although the pursued objectives were the same that the ones in this framework, the welfare of users and the cost of schedules, they were studied in a lexicographic way prioritizing the welfare of users over the cost of schedules. To solve the problem, a combination of the Adaptive Large Neighborhood Search (ALNS) and a custom heuristic method was used. 

In contrast to previous ones, in this work the  problem of the company is studied as a biobjective HCSP, that is, the goal is to study the trade off between the two objectives. Thus, the purpose is to obtain the Pareto frontier of the problem or a good approximation, which results in a set of solutions with different objective values from which the company can choose the one that best suits its needs. 

Different methods have been used in the literature to solve multiobjective HCSP. Some frameworks of multiobjective HCSP apply exact methods to solve small instances, such as the $\epsilon$-constraint approach \cite{Laumanns2006}, which consists of expressing some objectives as constraints. However, the complexity of this kind of problems (NP-Hard problems, since they are variants of VRP), requires the application of other techniques in realistic scenarios, such as metaheuristic/heuristic methods. 

Hence, to solve all the instances (both simulated and real instances), in this framework we design an algorithm based on the Multi-Directional Local Search (MDLS) metaheuristic method \cite{TRICOIRE2012}, which is based on exploring neighborhoods using single-objective searches. In addition, we use other techniques to show the good performance of this algorithm. On one hand, we develop an algorithm based on an improvement of the $\epsilon$-constraint approach and apply it to solve small instances. On the other hand, we adapt another metaheuristic technique, the Non-dominated Sorting Genetic Algorithm (NSGA-II) implemented in jMetal framework, to compare it with the MDLS algorithm. 

Section 2 is devoted to the mathematical model. In Section 3 the MDLS algorithm is explained in detail. Section 4 describes an improved version of the  
$\epsilon$-constraint method, the AUGMECON2, and the adaptation of the NSGA-II algorithm to this context. Finally, Section 5 deals with  computational experiments.

\section{The biobjective HCSP model} 

The goal of the problem under study is to correctly attend users needs, which should be carried out by the caregivers of the company. Each user requires a set of services with a predetermined duration, which depends on the corresponding task. Every service is assigned a time slot (morning, noon, afternoon or evening) on a specific day of the week. In this way, the user is allowed to indicate two types of time windows for any service: a hard time window establishing the time interval within which the service must be scheduled, and a soft time window, defining the time preferences of the user.

\begin{figure}[H]
	\begin{subfigure}{0.45\textwidth}
		\centering
		\includegraphics[scale=0.55,trim={1cm 0.9cm 1cm 0.8cm}, clip]{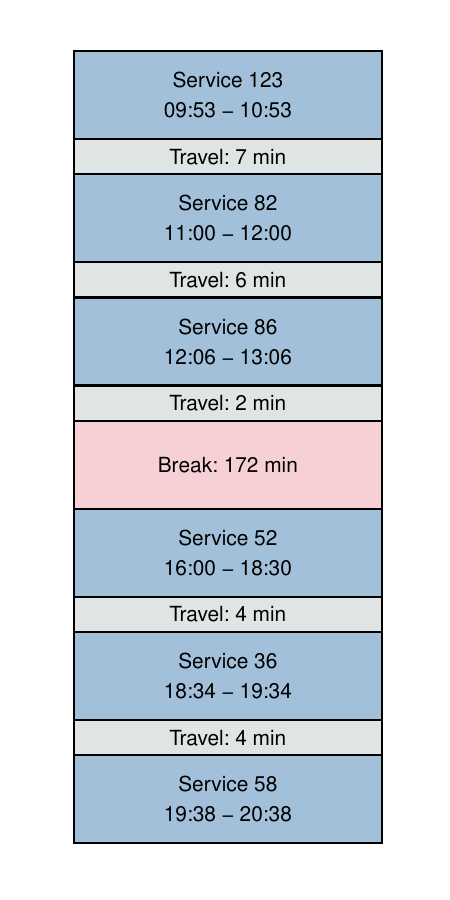}
		\caption{Schedule}
		\label{fig:schedule}
	\end{subfigure}
	\begin{subfigure}{0.5\textwidth}
		\centering
		\includegraphics[scale=0.9,trim={5cm 9cm 5cm 2cm}, clip]{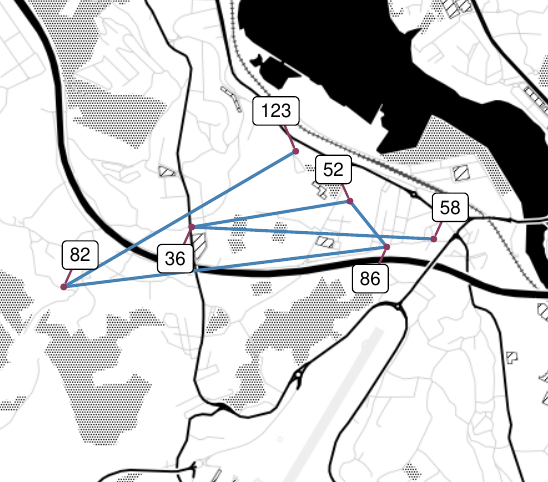}
		\caption{Route}
		\label{fig:route}
	\end{subfigure}
	\caption{Example of a working day for a caregiver.}
	\label{fig:caregiversday}
\end{figure}

Figure~\ref{fig:caregiversday} presents a real schedule and the associated route of a caregiver of the company in a working day, with 6 services of different durations and travel times between them.

The daily working time of the caregivers computes from the first service of the day to the last one. As it can be seen in Figure~\ref{fig:route}, the route of the caregiver starts at service 123 and ends at service 58. According to the company policy, the largest break of the day (i.e., the time they are not attending users nor travelling between their homes) should be discounted from the working day of the caregiver, as long as it is greater than or equal to 2 hours. The caregiver has a break between services 86 and 52 with a duration of 172 minutes (see Figure~\ref{fig:schedule}), which means that it will not be paid.

The maximum amount of time they can work each day, as well as the weekly agreed working time, are specified in the contract of each caregiver. The weekly time can be surpassed as long as the daily one is not exceeded and the caregivers are paid for their overtime. According to Figure~\ref{fig:schedule}, the caregiver has a scheduled working time of 7 hours and 43 minutes, which will not surpass the daily maximum of 8 hours.

With the aim of maintaining users' satisfaction, six levels of affinity are considered, determining how suitable is a caregiver to attend a user. These levels go from 0 to 5, where the higher the level, the more compatible is the assignation caregiver-user.

Next, we will formally introduce the sets, parameters and variables for the mathematical formulation of the optimization problem. 

\begin{table}[H]
	\renewcommand{\arraystretch}{0.95}
	\centering
	\begin{tabular}{ll}
		\hline & \\[-2.1ex]
		\textbf{Sets}    &     \\
		\hline & \\[-2.1ex]
		$D = \{1,...,7\}$      & Set of days.          \\
		$N=\{1,...,n\}$      & Set of caregivers.      \\
		$S=\{1,...,s-1\}$      & Set of services.      \\
		$S^{0} = S \cup \{0\}$      & Set of services and the initial dummy.      \\
		$S^{1} = S \cup \{s\}$      & Set of services and the ending dummy.      \\
		$S^{01} = S \cup \{0,s\}$      & Set of services and the initial and ending dummies.      \\
		$S_{-k} = S \setminus \{ k \}$      & Set of services except $k \in S$. Analogously, $S^{0}_{-k} = S^{0} \setminus \{ k \}$ and $S^{1}_{-k} = S^{1} \setminus \{ k \}$.       \\[4pt]
		\hline & \\[-2.1ex]
		\textbf{Parameters}    &     \\
		\hline & \\[-2.1ex]
		$\prm{\rho}{j}$      & It indicates if caregiver $i \in N$ can perform service $j \in S$.      \\
		$\prm{\lambda}{j}$      & Affinity level between caregiver $i \in N$ and service $j \in S$. \\
		$\eta_{j}$      & Duration of service $j \in S$.      \\
		$ [ \tw{\ubar{\alpha}}{j}, \tw{\bar{\alpha}}{j} ]$     & Hard time window of service $j \in S$ in day $d \in D$. \\
		& Note that, if service $j \in S$ does not belong to day $d \in D$, we set $ \tw{\ubar{\alpha}}{j} = \tw{\bar{\alpha}}{j} $.       \\
		$ [ \tw{\ubar{\beta}}{j}, \tw{\bar{\beta}}{j} ]$      & Soft time window of service $j \in S$ in day $d \in D$.      \\
		$ [ \var{\ubar{\gamma}}{}, \var{\bar{\gamma}}{}]$      & Availability time period of caregiver $i \in N$ at day $d \in D$.     \\
		$\theta_{jk}$      & Travel time between services $j \in S$ and $k \in S^{1}$.\\
		& Note that, if $k = s$, then $\theta_{js} = 0$.       \\
		$\prm{\nu}{}$      & Agreed weekly working time of caregiver $i \in N$. \\
		$\var{\nu}{}$ &  Maximum time caregiver $i \in N$ is allowed to work at day $d \in D$.      \\
		$\pi_{min}$      &  Minimum length of time required for the largest break to be unpaid.      \\[4pt]
		
		\hline
	\end{tabular}
\end{table}

\begin{table}[H]
	\renewcommand{\arraystretch}{0.95}
	\centering
	\begin{tabular}{ll}
		\hline & \\[-2.1ex]
		\textbf{Variables}    &     \\
		\hline & \\[-2.1ex]
		$\var{x}{jk}$      & It indicates if caregiver $i\in N$ goes from service $j\in S^{0}$ \\
		& to service $k\in S^1$ at day $d\in D$.      \\
		$\var{t}{j}$      &  For caregiver $i \in N$, it represents the starting time of service $j\in S^{01}$ \\
		& at day $d\in D$.       \\
		$\var{y}{jk}$      &  For caregiver $i\in N$, it indicates if the break between services $j \in S$ and $k \in S$  \\
		& has been selected to be discounted from the working day $d \in D$ .      \\
		$\var{\bar{y}}{}$      &  It states if there is no break for caregiver $i\in N$ at day $d\in D$.      \\
		$\var{r}{}$      & Greatest break of caregiver $i \in N$ at day $d  \in D$.\\
		$\var{u}{} $      & It indicates if the largest break of caregiver $i\in N$ at day $d\in D$ \\ 
		& is at least $\pi_{min}$.   \\
		$\var{\hat{r}}{}$      & Greatest break of caregiver $i\in N$ at day $d\in D$ if it is at least $\pi_{min}$.\\
		&  Otherwise, it will be 0.        \\
		$\prm{z}{} $      & Amount of overtime of caregiver $i \in N$.      \\
		$v_{j}^{start}$      & Penalization for carrying out service $j \in S$ before its soft time window.      \\
		$v_{j}^{end}$      & Penalization for carrying out service $j \in S$ after its soft time window.      \\
		\hline
	\end{tabular}
\end{table}

The objective functions and the constraints of the problem are:

{\allowdisplaybreaks
	\begin{alignat}{2}
	{f_1 }  = &
	 \textrm{ min } {\omega_1 \sum_{i \in N } \prm{z}{}} + \omega_2 \sum_{i \in N } \sum_{d \in D} (\var{t}{s} - \var{t}{0} - \var{\hat{r}}{}) \label{objective1} \\
	{f_2 }  =  &
	\textrm{ min } \omega_3 \sum_{i \in N } \sum_{d \in D } \sum_{j \in S } \sum_{k \in S^{1} }  \prm{\lambda}{j} \var{x}{jk}  + \omega_4 \sum_{j \in S } (v_{j}^{start} + v_{j}^{end}) \label{objective2}
	\end{alignat}}

Subject to
{\allowdisplaybreaks
	\begin{alignat}{2}
	\label{constraint05}
	&\sum_{i\in N} \sum_{d\in D} \sum_{k \in S^{1}_{-j}} \var{x}{jk} = 1 & \hspace*{2cm} & \forall j \in S \\
	\label{constraint06}
	&\sum_{i\in N} \sum_{d\in D} \sum_{j \in S^{0}_{-k}} \var{x}{jk} = 1 & \hspace*{2cm} & \forall k \in S \\
	\label{constraint07}
	&\sum_{d\in D} \sum_{k \in S^{1}_{-j}} \var{x}{jk} \le \prm{\rho}{j} & \hspace*{2cm} & \forall i\in N, \forall j \in S \\
	\label{constraint08}
	&\sum_{k \in S^{1}} \var{x}{0k} = 1 & \hspace*{2cm} & \forall i\in N, \forall d\in D \\
	\label{constraint09}
	& \sum_{j \in S^{0}} \var{x}{js} = 1 & \hspace*{2cm} & \forall i\in N, \forall d\in D  \\
	\label{constraint10}
	& \sum_{j \in S^{0}_{-h}} \var{x}{jh} - \sum_{k \in S^{1}_{-h}} \var{x}{hk}= 0 & \hspace*{2cm} & \forall i\in N, \forall d\in D, \forall h\in S \\
	\label{constraint11}
	&\tw{\ubar{\alpha}}{j}\sum_{k \in S^{1}_{-j}} \var{x}{jk}\le \var{t}{j} \le ( \tw{\bar{\alpha}}{j} - \eta_{j})\sum_{k \in S^{1}_{-j}} \var{x}{jk} & \hspace*{2cm} & \forall i\in N, \forall d\in D, \forall j\in S \\
	\label{constraint12}
	& \var{t}{j}+(\eta_{j} + \theta_{jk})\var{x}{jk}\le \var{t}{k}+\tw{\bar{\alpha}}{j}(1-\var{x}{jk}) & \hspace*{2cm} & \forall i\in N, \forall d\in D, \\
	\nonumber & & & \forall j\in S, \forall k\in S^{1}, j\neq k \\
	\label{constraint13}
	& \var{t}{0} \ge \var{\ubar{\gamma}}{} & \hspace*{2cm} & \forall i\in N, \forall d\in D \\
	\label{constraint14}
	& \var{t}{s} \le \var{\bar{\gamma}}{} & \hspace*{2cm} & \forall i\in N, \forall d\in D \\
	\label{constraint15}
	& \var{t}{0}\le \var{t}{k}+\var{\bar{\gamma}}{}(1-\var{x}{0k}) & \hspace*{2cm} & \forall i\in N, \forall d\in D,  \forall k\in S^{1} \\
	\label{constraint16}
	& \var{t}{0}\ge \var{t}{k}-\var{\bar{\gamma}}{}(1-\var{x}{0k}) & \hspace*{2cm} & \forall i\in N, \forall d\in D,  \forall k\in S^{1} \\
	\label{constraint17}
	& \var{t}{s}\le (\var{t}{j} + \eta_{j})+\var{\bar{\gamma}}{}(1-\var{x}{js}) & \hspace*{2cm} & \forall i\in N, \forall d\in D,  \forall j\in S \\
	\label{constraint18}
	& \var{t}{s}-\var{t}{0} - \var{\hat{r}}{} \le \var{\nu}{}  & \hspace*{2cm} & \forall i\in N, \forall d\in D \\
	%
	\label{constraint20}
	&\prm{z}{} \ge \sum_{d\in D} \left(\var{t}{s}-\var{t}{0} - \var{\hat{r}}{}\right) - \prm{\nu}{con} & \hspace*{2cm} & \forall i\in N \\
	\label{constraint21}
	&\var{r}{} \ge \var{t}{k}-(\var{t}{j} + \eta_{j} + \theta_{jk}) - \var{\bar{\gamma}}{}(1-\var{x}{jk}) & \hspace*{2cm} & \forall i\in N, \forall d\in D, \\
		\nonumber & & &   \forall j\in S, \forall k\in S, j\neq k \\
	\label{constraint22}
	&\var{r}{} \le \var{t}{k}-(\var{t}{j} + \eta_{j} + \theta_{jk}) +\var{\bar{\gamma}}{}(1-\var{x}{jk}) +\var{\bar{\gamma}}{}(1-\var{y}{jk}) & \hspace*{1.5cm} & \forall i\in N, \forall d\in D, \\
	\nonumber & & & \forall j\in S, \forall k\in S, j\neq k \\
	\label{constraint23}
	&\var{r}{} \le \var{\bar{\gamma}}{}(1-\var{\bar{y}}{}) & \hspace*{1.5cm} & \forall i\in N, \forall d\in D \\
	\label{constraint24}
	&\sum_{j\in S}\sum_{k\in S_{-j}} \var{y}{jk} + \var{\bar{y}}{}=1 & \hspace*{2cm} & \forall i\in N, \forall d\in D \\
	\label{constraint25}
	& \var{y}{jk} \leq \var{x}{jk} & \hspace*{2cm} & \forall i\in N, \forall d\in D, \\
		\nonumber & & &  \forall j\in S, \forall k\in S, j\neq k \\
	\label{constraint26}
	&\var{r}{} - \pi_{min} \geq \pi_{min} (\var{u}{}-1) & \hspace*{2cm} & \forall i\in N, \forall d\in D \\
	\label{constraint27}
	&\var{r}{} - \pi_{min} + \varepsilon \leq (\var{\bar{\gamma}}{}-\var{\ubar{\gamma}}{}) \var{u}{} & \hspace*{2cm} & \forall i\in N, \forall d\in D \\
	\label{constraint28}
	&\var{\hat{r}}{} \leq (\var{\bar{\gamma}}{}-\var{\ubar{\gamma}}{})\var{u}{} & \hspace*{2cm} & \forall i\in N, \forall d\in D \\
	\label{constraint29}
	&\var{\hat{r}}{} \leq \var{r}{} & \hspace*{2cm} & \forall i\in N, \forall d\in D \\
	\label{constraint30}
	&\var{\hat{r}}{} \geq \var{r}{} - (\var{\bar{\gamma}}{}-\var{\ubar{\gamma}}{}) (1- \var{u}{})& \hspace*{2cm} & \forall i\in N, \forall d\in D \\
	\label{constraint31}
	&v_{j}^{start} \geq \sum_{d \in D } \left(\tw{\ubar{\beta}}{j} \sum_{i \in N }\sum_{k \in S^{1}_{-j}}\var{x}{jk} - \sum_{i \in N }\var{t}{j} \right) & \hspace*{2cm} & \forall j\in S \\
	\label{constraint32}
	&v_{j}^{end} \geq \sum_{d \in D } \left(\sum_{i \in N }\var{t}{j} + (\eta_{j} - \tw{\bar{\beta}}{j})\sum_{i \in N }\sum_{k \in S^{1}_{-j}}\var{x}{jk} \right) & \hspace*{2cm} & \forall j\in S \\
	\label{constraint33}
	&\var{x}{jk} \in \{0,1\} & \hspace*{2cm} & \forall i\in N, \forall d\in D, \\
		\nonumber & & &  \forall j\in S^{0}, \forall k\in S^{1}, j\neq k\\
	\label{constraint34}
	&\var{y}{jk}\in \{0,1\} & \hspace*{2cm} & \forall i\in N, \forall d\in D, \\
	\nonumber & & &  \forall j\in S, \forall k\in S, j\neq k\\
	\label{constraint35}
	&\var{\bar{y}}{}, \var{u}{} \in \{0,1\}; \; \var{r}{}, \var{\hat{r}}{} \ge 0 & \hspace*{2cm} & \forall i\in N, \forall d\in D\\
	\label{constraint36}
	& \prm{z}{} \ge 0  & \hspace*{2cm} & \forall i\in N\\
	\label{constraint37}
	&\var{t}{j} \ge 0 & \hspace*{2cm} & \forall i\in N, \forall j\in S^{01}, \forall d\in D \\
	\label{constraint38}
	& \var{r}{}, \var{\hat{r}}{} \ge 0 & \hspace*{2cm} & \forall i\in N, \forall d\in D \\
	\label{constraint39}
	& v_{j}^{start}, v_{j}^{end} \ge 0 & \hspace*{2cm} &  \forall j\in S
	\end{alignat}}

There are two objectives in this problem: the cost of schedules and the welfare of users. 
The cost of schedules is addressed in Objective~(\ref{objective1}) and it is divided into two elements: the overtime of caregivers and their worked time. Since both of them depend on the same time units, we consider their weights as $\omega_1 = \omega_2 = 1$.
On the other hand, Objective~(\ref{objective2}) accounts for the welfare and includes the penalization of carrying out the services outside their soft times and the affinity caregivers-services. Since we want to prioritize the affinity over the soft time windows penalization we set $ \omega_3 = - \max \{ 1 ,  \sum_{j \in S} [ \max_{d \in D} \{ \tw{\ubar{\beta}}{j} - \tw{\ubar{\alpha}}{j}\} + \max_{d \in D} \{ \tw{\bar{\alpha}}{j} - \tw{\bar{\beta}}{j}\} ] \} $ and  $\omega_4  = 1$\footnote{This is based on the fact that $v_{j}^{start} \leq \max_{d \in D} \{ \tw{\ubar{\beta}}{j} - \tw{\ubar{\alpha}}{j}\}$ and $v_{j}^{end} \leq \max_{d \in D} \{ \tw{\bar{\alpha}}{j} - \tw{\bar{\beta}}{j}\}$ for all $j \in S$.}.

Regarding the constraints, (\ref{constraint05}) - (\ref{constraint06})  guarantee that all services will be assigned to an available caregiver. Constraints (\ref{constraint08}) - (\ref{constraint09}) state that every caregiver starts and finishes each working day with the dummy services. Constraint (\ref{constraint10}) checks that routes are not segmented. Constraint (\ref{constraint11}) imposes that services are performed within their required time windows. Travel times and caregivers time windows availability are upheld in Constraints (\ref{constraint12}) - (\ref{constraint14}). The schedule of the dummy services is controlled in Constraints (\ref{constraint15}) - (\ref{constraint17}). Constraint (\ref{constraint18}) ensures that the daily maximum number of working hours will not be exceeded. The lower bound of the overtime for every caregiver is set at Constraint (\ref{constraint20}).
Constraints (\ref{constraint21}) - (\ref{constraint25}) are related with the computation of the single daily largest break for each caregiver. Constraints (\ref{constraint26}) - (\ref{constraint27}) check whether the largest daily break of the caregiver should be remunerated or not. In case this break is not paid, Constraints (\ref{constraint28}) - (\ref{constraint30}) claim that it will be discounted from working time.
The time to carry a service outside its soft time window is obtained in Constraints (\ref{constraint31}) - (\ref{constraint32}). 
Lastly, Constraints (\ref{constraint33})-(\ref{constraint39}) establish variables domain.

\section{The main approach: the BIALNS algorithm}
In this section we propose an algorithm to solve the biobjective HSCP of Section 2. The fundamentals of the algorithm are based on the Multi-Directional Local Search (MDLS) metaheuristic method \cite{TRICOIRE2012}. In this context, it will be necessary to obtain solutions that prioritize the welfare over the cost, and vice versa. To this aim, an Adaptive Large Neighbourhood Search (ALNS) algorithm is combined with a customized scheduling heuristic method.

Algorithm~\ref{alg:Biobj} shows the method designed to solve the biobjective problem. A solution of the problem is $\pmb{\omega} = (\pmb{x},\pmb{t})$ and it is composed by the routes ($\pmb{x}$) and the schedules ($\pmb{t}$). The routes describe the order of services to be performed by each caregiver and the schedules are the starting times of each service. The goal of the algorithm is to obtain the set of non dominated solutions ($\Omega$), so it relies on a set of solutions that are composed by different routes ($\hat{\Omega}$). 

\vspace*{0.5cm}
\begin{algorithm}[H]
	\setstretch{1.05}
	\SetKwComment{Comment}{// }{}
	\KwData{Services ($S$), Caregivers ($N$), Lexicographic function welfare-cost ($f_{wc}$), Lexicographic function cost-welfare ($f_{cw}$)}
	\Comment{Initialize the non dominated set}
	$\pmb{\omega}_{wc}$ $\gets$ \textit{{\textbf{initialSolution}}}($S$, $N$, $f_{wc}$) \Comment{Get initial solution for welfare-cost}
	$\pmb{\omega}_{cw}$ $\gets$ \textit{{\textbf{initialSolution}}}($S$, $N$, $f_{cw}$) \Comment{Get initial solution for cost-welfare}
	$\pmb{\omega}_{wc}, \hat{\Omega}$ $\gets$ \textit{{\textbf{ALNS}}}($f_{wc}$, $\pmb{\omega}_{wc}$, $\hat{\Omega} = \emptyset$) \Comment{ALNS for the lexicographic welfare-cost}
	$\pmb{\omega}_{cw}, \hat{\Omega}$ $\gets$ \textit{{\textbf{ALNS}}}($f_{cw}$, $\pmb{\omega}_{cw}$, $\hat{\Omega}$) \Comment{ALNS for the lexicographic cost-welfare}
	$\Omega$ $\gets$ \textit{{\textbf{updateNonDominatedSet}}}($\Omega$, $\pmb{\omega}_{wc}$) \Comment{Get non dominated solutions}
	$\Omega$ $\gets$ \textit{{\textbf{updateNonDominatedSet}}}($\Omega$, $\pmb{\omega}_{cw}$) \Comment{Get non dominated solutions}
	
	\Comment{Obtain solutions composed by different routes}
	\While{stopping criteria not met}{
		$\bar{\pmb{\omega}}$ $\gets$ \textit{{\textbf{chooseRandomSolution}}($\hat{\Omega}$, $\Omega$)} \Comment{Get solution to modify}
		$\pmb{\omega}_{wc}, \hat{\Omega}$ $\gets$ \textit{\hyperref[sec:ALNSBI]{\textbf{ALNS}}}($S$, $N$, $f_{wc}$, $\bar{\pmb{\omega}}$) \Comment{ALNS for the lexicographic welfare-cost}
		$\pmb{\omega}_{cw}, \hat{\Omega}$ $\gets$ \textit{\hyperref[sec:ALNSBI]{\textbf{ALNS}}}($S$, $N$, $f_{cw}$, $\bar{\pmb{\omega}}$) \Comment{ALNS for the lexicographic cost-welfare}
		$\Omega$ $\gets$ \textit{\hyperref[sec:updateNonDominatedSet]{\textbf{updateNonDominatedSet}}}($\Omega$, $\pmb{\omega}_{wc}$) \Comment{Update the non dominated solutions}
		$\Omega$ $\gets$ \textit{\hyperref[sec:updateNonDominatedSet]{\textbf{updateNonDominatedSet}}}($\Omega$, $\pmb{\omega}_{cw}$) \Comment{Update the non dominated solutions}
	}
	\Comment{Get non dominated solutions}
	\While{stopping criteria not met}{
		$\bar{\pmb{\omega}}$ $\gets$ \textit{{\textbf{chooseRandomSolution}}($\hat{\Omega}$, $\Omega$)} \Comment{Get solution to modify}
		\textit{$\Omega$} $\gets$ \textit{{\textbf{modifyScheduleWelf}}($\bar{\pmb{\omega}}, \Omega$)} \Comment{Modify the schedule to improve welfare}
		\textit{$\Omega$} $\gets$ \textit{{\textbf{modifyScheduleCost}}($\bar{\pmb{\omega}}, \Omega$)} \Comment{Modify the schedule to improve cost}
	}
	
	\Return $\Omega$
	\caption{\textbf{BIALNS algorithm}  } \label{alg:Biobj}
\end{algorithm}
\vspace*{0.5cm}

The algorithm is divided into three steps. In the first one,  the initial solutions (lines 1 - 2) are obtained by applying a  random greedy insertion operator, either prioritizing welfare over cost or cost over welfare. To enhance the solutions, the ALNS philosophy adapted to our particular context (Algorithm~\ref{alg:ALNSBI}) is implemented for each lexicographic objective (lines 3 - 4). Then, these solutions are used to update the set of non dominated solutions (lines 5 - 6).

The second step consists in generating solutions composed by different routes. Thus, a solution is selected at random from the ones already found (line 8) and then it is updated using the ALNS method for both lexicographic objectives (lines 9 - 10). The set of solutions composed by different routes is updated during the ALNS. The new solutions are used to update the set of non dominated solutions (lines 11 - 12).

Finally, during the third step, non dominated solutions are generated. To this aim, a solution is chosen at random at random (line 14), modifying its schedule in order to improve welfare (line 15) and cost (line  16). The methods used to modify the schedule are described in Subsections~\ref{sec:imprwelf}~and~\ref{sec:imprcost}. 

\subsection{ALNS method}

Algorithm~\ref{alg:ALNSBI} describes the ALNS method. The input data are the initial solution ($\omega$), the removal ($\Sigma_{rem}$) and insertion ($\Sigma_{ins}$) operators, their weights ($\sigma_{rem}$, $\sigma_{ins}$) and the set of solutions composed by different routes ($\hat{\Omega}$).

\vspace*{0.5cm}
\begin{algorithm}[H]
	\setstretch{1.05}
	\SetKwComment{Comment}{// }{}
	\KwData{Objective function ($f$), initial solution ($\omega$), set of solutions ($\hat{\Omega}$)}
	\caption{\textbf{ALNS -} Adaptive Large Neighbourhood search} \label{alg:ALNSBI} 
	\textit{$\sigma_{rem}$} $\gets$ \textit{$(1,...,1)$},  
	\textit{$\sigma_{ins}$} $\gets$ \textit{$(1,...,1)$}, $\pmb{\omega}'$ $\gets$ $\pmb{\omega}$ \\
	\Comment{Improve the solution}
	\While{stopping criteria not met}{
		$\varsigma_{rem}$ $\gets$ \textit{{\textbf{chooseRandom}}}($\sigma_{rem}, \Sigma_{rem}$), $\varsigma_{ins}$ $\gets$ \textit{{\textbf{chooseRandom}}}($\sigma_{ins}, \Sigma_{ins}$) \Comment{Get operators}
		$\bar{\pmb{\omega}}$ $\gets$ \textit{{\textbf{destroySolution}}}($\pmb{\omega}$, $\varsigma_{rem}$, $f$), $\pmb{\omega}^*$ $\gets$ \textit{{\textbf{repairSolution}}}($\bar{\pmb{\omega}}$, $\varsigma_{ins}$, $f$) \Comment{Obtain new solution}
		\If{f($\pmb{\omega}^*$) $<$ f($\pmb{\omega}'$)}{
			\textit{$\pmb{\omega}'$} $\gets$ \textit{$\pmb{\omega}^*$} \Comment{Update best solution}
		}
		\Comment{Update the set of multiple routes}
		$add$ $\gets$ $true$\\
		\For{$\hat{\pmb{\omega}} \in \hat{\Omega}$}{
			\If(\tcp*[h]{Check if the routes are equal}){$\pmb{x} = \hat{\pmb{x}}$}{
				$add$ $\gets$ $false$\\
				\textbf{break} \\
			}
		}
		
		\If{$add = true$}{
			$\hat{\Omega}$ $\gets$ $\hat{\Omega} \cup \{\pmb{\omega}^*\}$ \Comment{Add the solution to the set}
		}
		\Comment{Update current solution}
		$\pmb{\omega}$ $\gets$ \textit{{\textbf{acceptanceCriteria}}}($\pmb{\omega}^*$, $\pmb{\omega}'$)\\
		\Comment{Update the weights of the operators}
		\textit{$\sigma_{rem}$} $\gets$ \textit{{\textbf{updateWeights}}}($\sigma_{rem}$, $f$, $\pmb{\omega}'$, $\pmb{\omega}^*$),
		\textit{$\sigma_{ins}$} $\gets$ \textit{{\textbf{updateWeights}}}($\sigma_{ins}$, $f$, $\pmb{\omega}'$, $\pmb{\omega}^*$)\\
	}

	\Return $\pmb{\omega}', \hat{\Omega}$
	
\end{algorithm}
\vspace*{0.5cm}

The ALNS consists of destroying and repairing a given solution using randomly chosen operators (lines 3 - 4). The newly generated solution is used to update the best solution found so far (lines 5 - 6). Then, we check if the routes of the new solution differ from the ones on the set, in which case the solution is added to the set (lines 7 - 13). Finally, the current solution and the weights of the operators are updated (lines 14 - 15).

The new solution is accepted with probability $exp(-(f(\pmb{\omega}^*) - f(\pmb{\omega}'))/T_{i})$, where $T_{i}$ is the temperature that decreases after each iteration $i$ according to formula $T_i=\beta^{i} T_{0}$, a cooling parameter $0 < \beta < 1$ and the initial temperature $T_0 > 0$. The weights of the removal and insertion operators are increased, after each iteration, if they generated an improved solution. Table~\ref{table:operators} describes all the operators used in the ALNS method.

Since a HCSP is a routing and scheduling problem, to evaluate the insertion operators it is not only necessary to know the routes of each caregiver, but also the starting times of every service. Because of this, after obtaining a new route using any of the insertion operators, it is necessary to define its schedule. 

\begin{table}[H]
	\renewcommand{\arraystretch}{1.1}
	\centering
	\begin{tabular}{ll}
		\hline
		\textbf{Operator}    &  \textbf{Description}    \\
		\hline
		Random removal     & Services to remove from the route are selected at random.          \\ \hline
		Related removal      & This operator iteratively removes from the route the service that is   \\ 
		& most related to a random already removed service.\\ \hline
		Cost removal      & Services that contribute the most to the objective function value  \\
		& are removed from the routes.      \\ \hline
		1-Route removal      & Randomly selected routes are removed from the solution \\
		& until the required number of services have been deleted.     \\
		\hline
		2-Route removal      & Two routes selected at random are removed from the solution.     \\
		\hline
		Basic greedy (BG) insertion     &  Service that results in the least objective function increase is \\
		& added to the solution.\\ \hline
		Random greedy (RG) insertion    & A randomly selected service is scheduled in the best possible position.\\ \hline
		Different caregiver BG     & Basic greedy insertion but trying to guarantee that services will \\
		insertion &  be assigned to a different caregiver.      \\ \hline
		Different caregiver RG    &  Basic greedy insertion but trying to guarantee that services will \\
		insertion &  be assigned to a different caregiver.      \\
		\hline
		
	\end{tabular}
	\caption{Removal and insertion operators.}
	\label{table:operators}
\end{table}

\begin{figure}[H]
	\begin{center}
		\includegraphics[width=1\textwidth,trim={1cm 11.4cm 4.7cm 1.0cm}, clip]{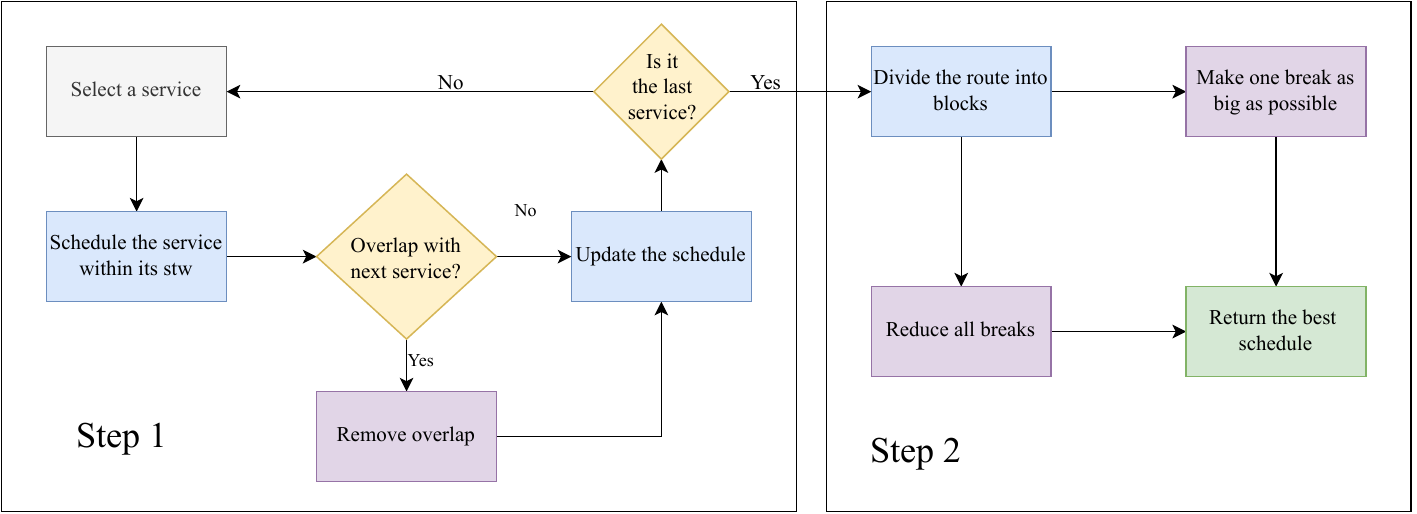}
	\end{center}
	\caption{Scheme to prioritize welfare over  cost.} \label{fig:diagrama_algwc}
\end{figure}

When prioritizing welfare over cost, the schedule of a route is obtained using the scheme of the  heuristic algorithm presented in Figure~\ref{fig:diagrama_algwc} (a complete description of this method can be found in \citet{MendezFernandez2023}). The algorithm is divided into two steps. In the first one,  schedule with best penalization value is found. In the second step, schedule is modified in order to improve the cost. 

\begin{figure}[H]
	\begin{center}
		\includegraphics[width=1\textwidth,trim={1cm 11.6cm 5.2cm 1.0cm}, clip]{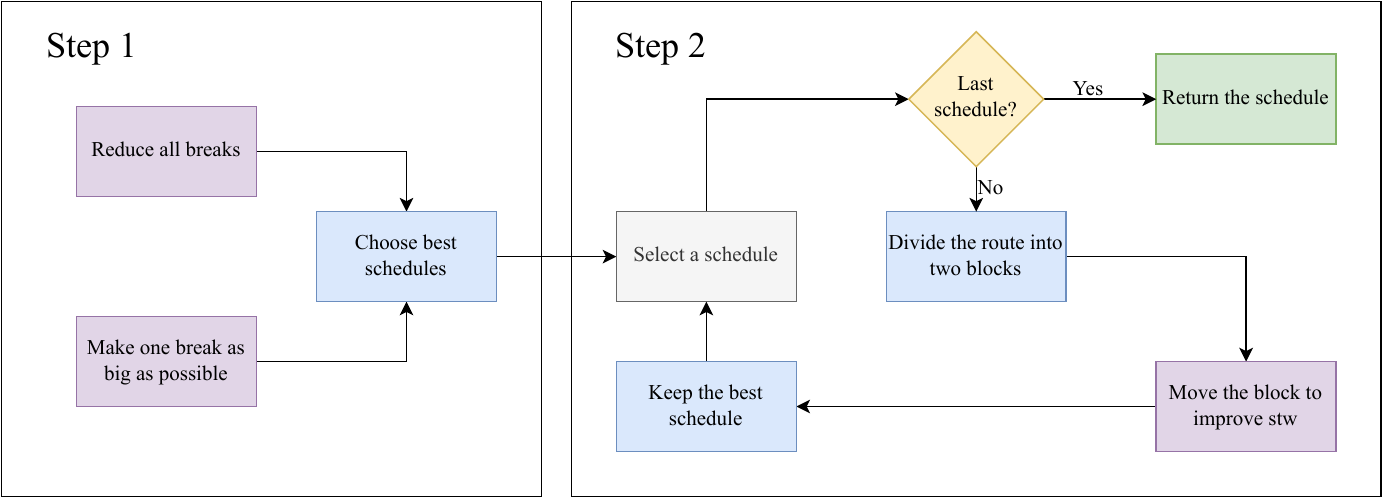}
	\end{center}
	\caption{Scheme to prioritize cost over  welfare.} \label{fig:diagrama_algcw}
\end{figure}

The scheme of the algorithm to obtain a route  schedule prioritizing cost over welfare is presented in Figure~\ref{fig:diagrama_algcw} (a full description of the method can be found in \citet{MendezFernandez_thesis}). 

The algorithm is divided into two steps: in the first one schedules with best cost value are found and, in the second step, schedules are modified in order to improve the soft time window penalization.

\subsection{Modify the schedule to improve the welfare}\label{sec:imprwelf}

\begin{figure}[H]
	\begin{center}
		\includegraphics[scale=0.75,trim={0cm 11.3cm 0cm 4.2cm}, clip]{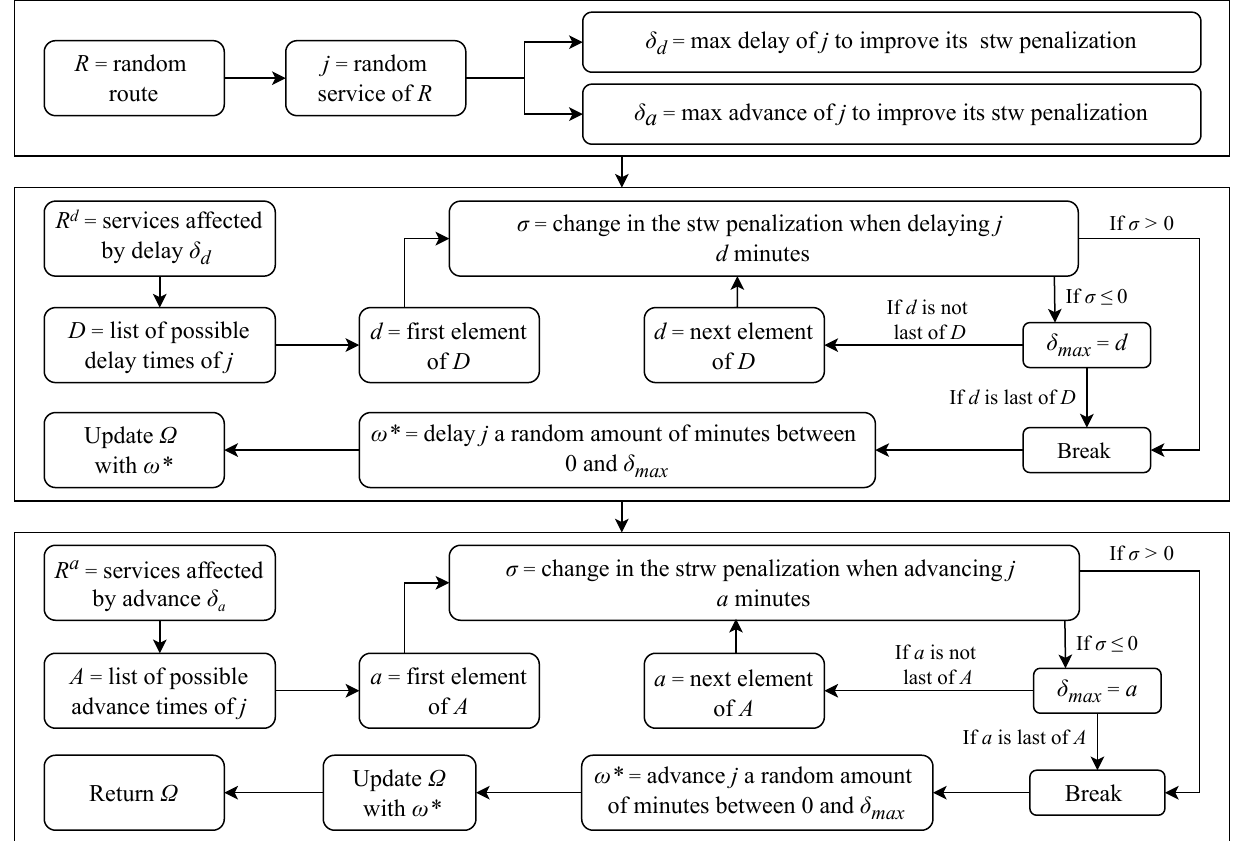}
	\end{center}
	\caption{Scheme to improve welfare.} \label{fig:improve_welf_alg}
\end{figure}

The algorithm developed to modify a schedule in order to improve the welfare, which only consists in the preferred time window penalization because the route is fixed, is shown in Figure~\ref{fig:improve_welf_alg} and thoroughly described in \citet{MendezFernandez_thesis}.

The algorithm is divided into three phases:
\begin{enumerate}
	\setlength\itemsep{-0.25em}
	\item In the first step a route and a service are randomly selected. Then maximum times to advance and delay the service, so that penalization is not increased, are computed.
	\item The maximum time to delay the service, without increasing the route  penalization, is obtained. The service is then delayed a random amount of time between 0 and the maximum. Then, non dominated set is updated.
	\item The maximum time to advance the service, without increasing the route  penalization, is obtained. The service is then advanced a random amount of time between 0 and the maximum. Finally, non dominated set is updated.
\end{enumerate}

Next, we propose an example to explain the practical application of the algorithm.

\begin{example} \label{ex:imprwelf}

	Figure~\ref{fig:move_welf0} shows an example of a modification of the schedule of a caregiver route to improve welfare. The route is composed of 6 services with its corresponding hard and soft time windows. The cost of the initial schedule is $f_1 = 570$ minutes and its soft time window penalization is $f_2= 30$ minutes\footnote{Note that in this case both objectives are measured in minutes because the route is already fixed. Moreover, once the route is fixed, the affinity is disregarded in the computation of $f_2$.}.
	
	\begin{figure}[h]
		\begin{center}
			\includegraphics[scale=0.65,trim={0cm 0cm 0cm 0cm}, clip]{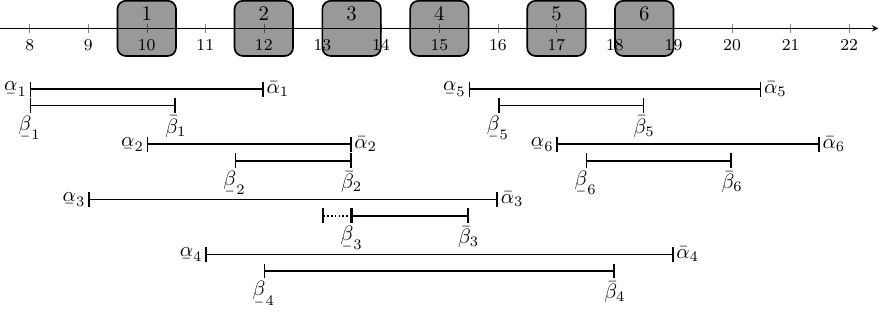}
		\end{center}
		\caption{Initial schedule.} \label{fig:move_welf0}
	\end{figure}
	
	The earliest and latest starting times of the services are: $t_1^e = 0$, $t_2^e = 120$, $t_3^e = 180$, $t_4^e = 240$, $t_5^e = 450$, $t_6^e = 540$ and $t_1^l =  180 $, $t_2^l = 270 $, $t_3^l = 450 $, $t_4^l = 600 $, $t_5^l =  720 $, $t_6^l = 780$. 
	
	Let us consider that, for instance, our purpose is to move service 4. Since its earliest and latest times are in the soft time window, the service will be advanced or delayed without increasing its penalization.
	The maximum delay time is $\delta_d= \min\{\tw{\bar{\beta}}{4} - \eta_4 - t_4, t_4^l - t_4\} = \min\{600 - 60 - 390, 600 - 390\} = 150$, while the maximum advance time is $\delta_a = \min\{t_4 - \tw{\ubar{\beta}}{4}, t_4 - t_4^e\} = \min\{390 - 240, 390 - 240\} = 150$.
	
	Next, the effect of the delay of service 4 in the schedule is studied. The delay times of service 4 that could result in a change of penalization are given by $D = \{120, 150\}$ (see Figure~\ref{fig:move_welf_posib}). To choose the best option it is necessary to compute the penalization. In Figure~\ref{fig:move_welf1} the penalization remains unchanged, because services 4, 5 and 6 are scheduled within their soft time windows. But, in the schedule presented in Figure~\ref{fig:move_welf2} the penalization of service 5 increases. Therefore, the best option is $\delta_{max} = 120$.
	
	\begin{figure}[H]
		\begin{subfigure}{0.5\textwidth}
			\centering
			\includegraphics[scale=0.55,trim={2.5cm 0cm 0cm 0cm}, clip]{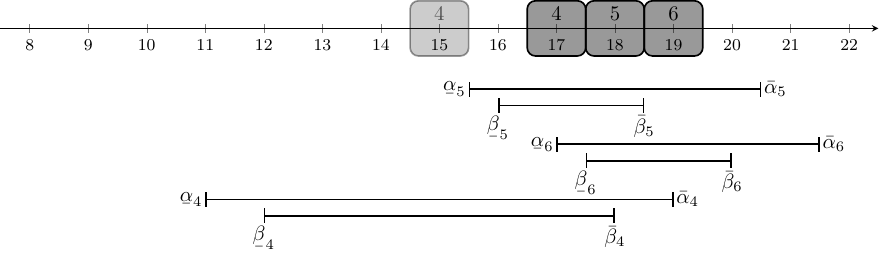}
			\caption{Delay of $120$ min.}
			\label{fig:move_welf1}
		\end{subfigure}
		\begin{subfigure}{0.5\textwidth}
			\centering
			\includegraphics[scale=0.55,trim={2.5cm 0cm 0cm 0cm}, clip]{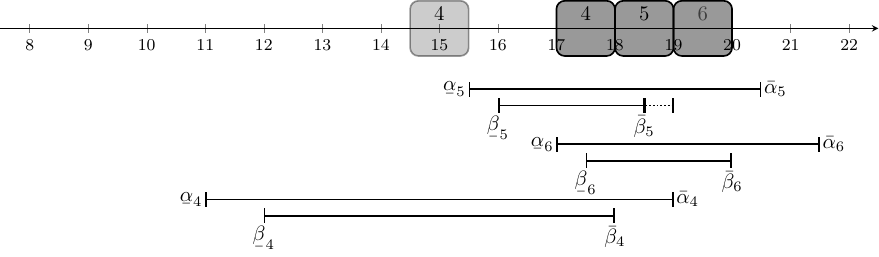}
			\caption{Delay of $150$ min.}
			\label{fig:move_welf2}
		\end{subfigure}
		\caption{Possible movements of service 4.}
		\label{fig:move_welf_posib}
	\end{figure}
	
	So, to delay the service a number between $\delta_{min} = 0$ and $\delta_{max} = 120$ is chosen, in this case $90$. The resulting schedule, in Figure~\ref{fig:move_welf3}, has a cost of $f^*_1 = 420$ and a soft time window penalization of $f^*_2=30$. That means that the new solution dominates the original one ($f^*_1 = 420 < f_1 = 570$ and $f^*_2=30 = f_2$).
	
	\begin{figure}[H]
		\begin{center}
			\includegraphics[scale=0.65,trim={0cm 0cm 0cm 0cm}, clip]{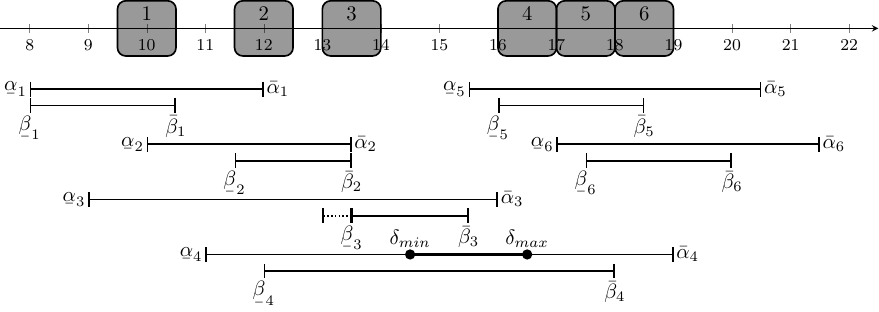}
		\end{center}
		\caption{Delayed schedule.} \label{fig:move_welf3}
	\end{figure}
	
	As far as the advance of service 4 is concerned, it could result in a change of penalization given by $A = \{ 150\}$ (see Figure~\ref{fig:move_welf4}). In this case, even though we can see that the penalization of services 2 and 3 increases, we set $\delta_{max} = 150$. Although the main objective is to improve the penalization, there is still a chance of increasing the penalization in case a non dominated solution is found.

	\begin{figure}[H]
		\begin{center}
			\includegraphics[scale=0.65,trim={0cm 0cm 2.8cm 0cm}, clip]{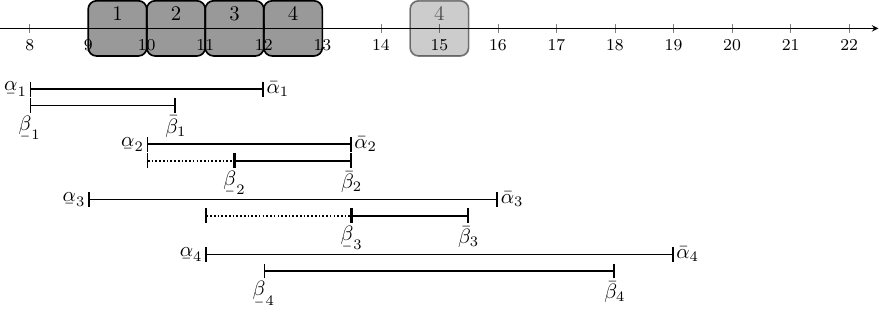}
		\end{center}
		\caption{Advance of $150$ min.} \label{fig:move_welf4}
	\end{figure}
	
	To advance the service we choose a number between $\delta_{min} = 0$ and $\delta_{max} = 150$, in this case $90$. The resulting schedule (see Figure~\ref{fig:move_welf5}) has a cost of $f'_1 = 450$ and a soft time window penalization of $f'_2=120$. This means that the new schedule is a non dominated solution ($f'_1 = 420 < f^*_1 = 450$ and $f'_2=120 > f^*_2=30$).
	
	\begin{figure}[H]
		\begin{center}
			\includegraphics[scale=0.65,trim={0cm 0cm 0cm 0cm}, clip]{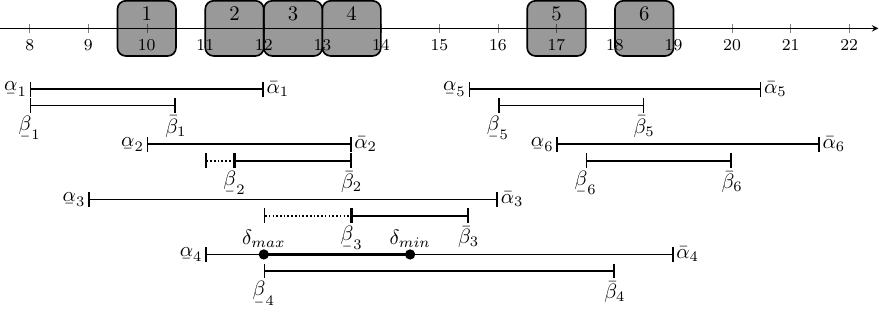}
		\end{center}
		\caption{Advanced schedule.} \label{fig:move_welf5}
	\end{figure}

\end{example}

\subsection{Modify the schedule to improve  cost}\label{sec:imprcost}

\begin{figure}[H]
	\begin{center}
		\includegraphics[scale=0.75,trim={0.0cm 5.0cm 0.0cm 1.0cm}, clip]{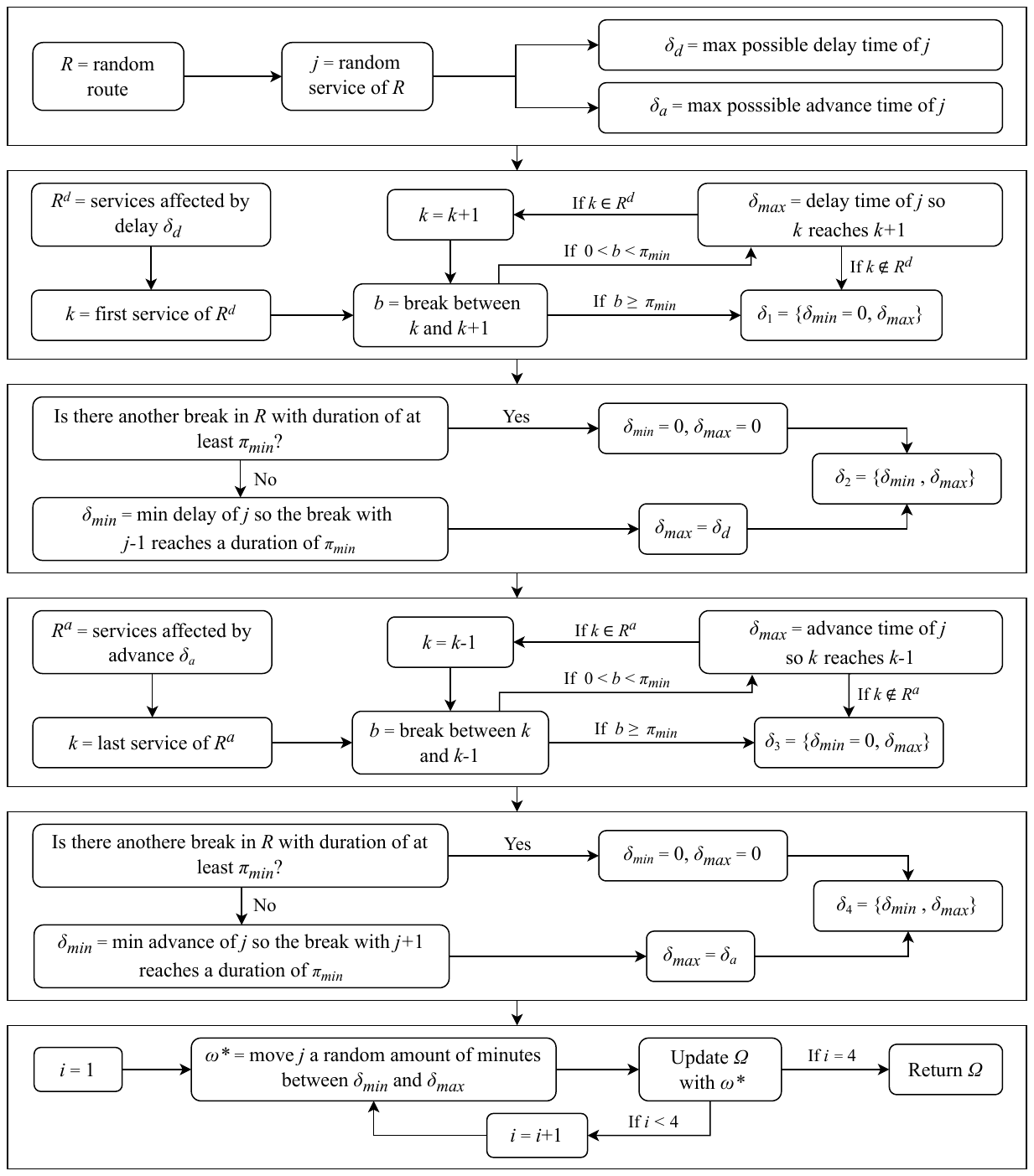}
	\end{center}
	\caption{Scheme to improve cost.} \label{fig:improve_cost_alg}
\end{figure}

The algorithm developed to modify a schedule in order to improve cost, which in this case only consists in modifying the breaks of the solution, is shown in Figure~\ref{fig:improve_cost_alg} and thoroughly described in \citet{MendezFernandez_thesis}.

The algorithm is divided into 6 steps:

\begin{itemize}
	\setlength\itemsep{-0.25em}
	\item The method starts by selecting a route and a service at random. After that, the maximum time to advance or delay a service is computed.
	\item The feasible delay of a service, in order to reduce the breaks that happen after it in the initial schedule of the route, is obtained.
	\item The maximum and minimum times to delay a service, to increase the duration of the break that happens right before it, are computed.
	\item The feasible advance of a service, in order to reduce the breaks that happen before it in the initial schedule of the route, is obtained.
	\item The maximum and minimum times to advance a service, to increase the duration of the break that happens right after it, are computed.
	\item Finally the service is moved a random amount of time, according to the times obtained before. That results in several solutions that are used to update the non dominated set.
\end{itemize}

Next example shows the modification of the schedule of Example~\ref{ex:imprwelf} to improve  cost.

\begin{example} \label{ex:impcost}
	After the modifications of the schedule suggested in Example~\ref{ex:imprwelf}, two non dominated solutions, with objective function values ($f^*_1 = 450$, $f^*_2=30$) and ($f'_1 = 420$, $f'_2=120$), are obtained\footnote{The original schedule was dominated.}.
	
	\begin{figure}[h]
		\begin{center}
			\includegraphics[scale=0.65,trim={0cm 0cm 0cm 0cm}, clip]{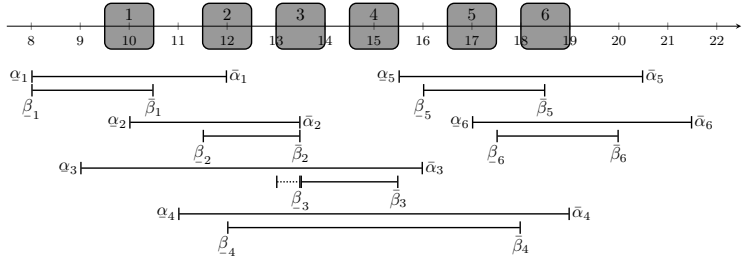}
		\end{center}
		\caption{Route that we want to modify to improve the cost.} \label{fig:move_cost0}
	\end{figure}
	
	The service to be moved is $j = 4$, its maximum delay time is $\delta_d = t^l_j - t_j = 600 - 390 = 210$ and the maximum advance time is $\delta_a = t_j - t^e_j = 390 - 240 = 150$ (line 12). The cost can be improved by moving service 4 according to the options described below.

	The first option is to delay the service, which would result in reducing the breaks after it. In this case, setting $\delta_{max} = 90$ completely removes the breaks after $j$ (lines 13 - 26). To delay the service, a number between $\delta_{max} = 0$ and $\delta_{max} = 90$ is chosen, in this case $60$. The resulting schedule, in Figure~\ref{fig:move_cost2}, has a cost of $f^1_1 = 570$ and a soft time window penalization of $f^1_2=30$. That means that the new solution is dominated ($f^1_1 = 570 > f^*_1 = 450 $ and $f^1_2= f^*_2=30 $).
	
	\begin{figure}[H]
		\begin{center}
			\includegraphics[scale=0.65,trim={0cm 0cm 0cm 0cm}, clip]{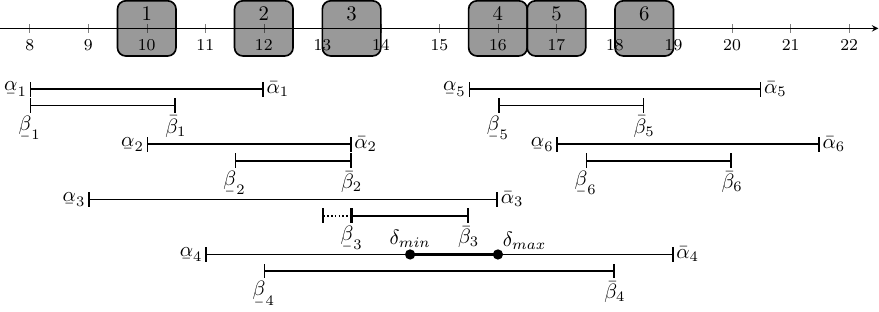}
		\end{center}
		\caption{Delayed schedule.} \label{fig:move_cost2}
	\end{figure}
	
	The second option is to delay the service in order to guarantee that the break before it will have a duration of $\pi_{max} = 120$ or more, which happens if we set $\delta_{min} = 90$ and $\delta_{max} = 210$. To delay the service, a number between $\delta_{max} = 90$ and $\delta_{max} = 210$ is chosen, in this case $180$. The resulting schedule, in Figure~\ref{fig:move_cost4}, has a cost of $f^2_1 = 450$ and a soft time window penalization of $f^2_2=150$. That means that the new solution is dominated ($f^2_1 = 450 = f^*_1$ and $f^2_2=150 > f^*_2 = 30$).
	
	\begin{figure}[H]
		\begin{center}
			\includegraphics[scale=0.65,trim={0cm 0cm 0cm 0cm}, clip]{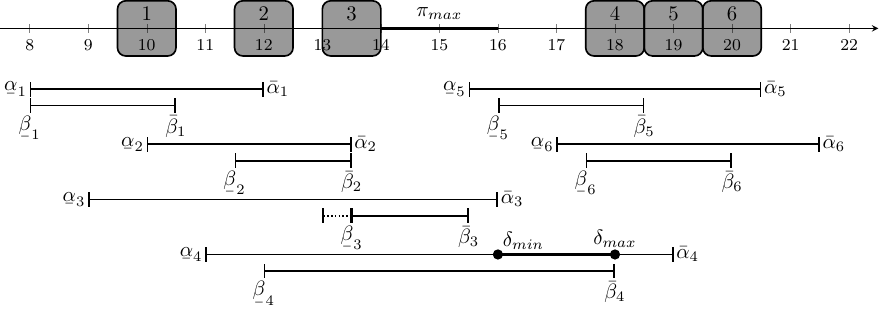}
		\end{center}
		\caption{Delayed schedule.} \label{fig:move_cost4}
	\end{figure}
	
	The third option consists of advancing the service to reduce the breaks that happen before it. If we set $\delta_{max} = 120$ we remove all breaks. To advance the service, a number between $\delta_{max} = 0$ and $\delta_{max} = 120$ is chosen, in this case $30$. The resulting schedule, Figure~\ref{fig:move_cost6}, has a cost of $f^3_1 = 570$ and a soft time window penalization of $f^3_2=30$. That means that the new schedule is dominated ($f^3_1 = 570 > f^*_1 = 450$ and $f^3_2= f^*_2=30 $).
	
	\begin{figure}[H]
		\begin{center}
			\includegraphics[scale=0.65,trim={0cm 0cm 0cm 0cm}, clip]{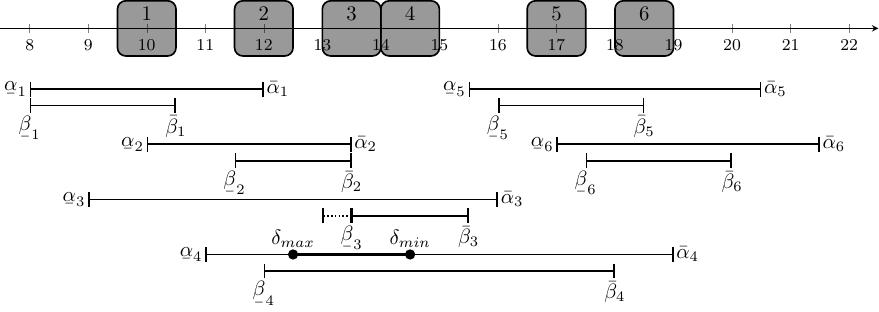}
		\end{center}
		\caption{Advanced schedule.} \label{fig:move_cost6}
	\end{figure}
	
	The final option is to get the advance times that make the break after the service of, at least, a duration of $\pi_{max} = 120$. It can be achieved by setting $\delta_{min} = 60$ and $\delta_{max} = 150$. To advance the service, a  number between $\delta_{min} = 60$ and $\delta_{max} = 150$ is chosen, in this case $120$. The resulting schedule, Figure~\ref{fig:move_cost8}, has a cost of $f^4_1 = 390$ and a soft time window penalization of $f^4_2=180$. That means that the new schedule is a non dominated solution ($f^4_1 = 390 < f^*_1 = 450$ and $f^4_2=180 > f^*_2=30$, $f^4_1 = 390 < f'_1 = 420$ and $f^4_2=180 > f'_2=120$).
	
	\begin{figure}[H]
		\begin{center}
			\includegraphics[scale=0.65,trim={0cm 0cm 0cm 0cm}, clip]{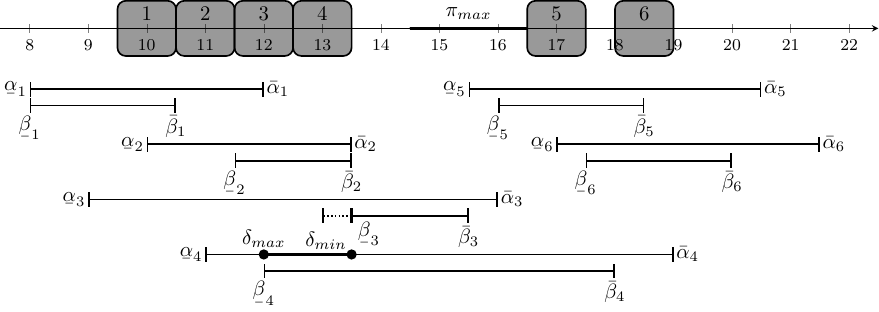}
		\end{center}
		\caption{Advanced schedule.} \label{fig:move_cost8}
	\end{figure}

\end{example}

\section{The AUGMECON2 method and the NSGA-II-based algorithm}
This section is devoted to introduce other two known state of the art multiobjective methods, with the objective of analysing the relative quality of the BIALNS algorithm. First, an exact method is described, the AUGMECON2, which is an improved version of the $\epsilon$-constraint method. Then, an alternative metaheuristic algorithm is proposed, based on the Non-Dominated Sorting Genetic Algorithm II (NSGA-II) technique. 

\newpage
\subsection{The AUGMECON2}
Now we describe the version of the AUGMECON2 method specific to our biobjective problem, $P$, which is:
\begin{alignat*}{2}
    \textrm{ min } &
    f_1(\pmb{x}), f_2(\pmb{x})  \\
    \textrm{ st. } & \pmb{x} \in F 
\end{alignat*}

Given two feasible solutions $\pmb{x}$ and $\pmb{y}$, it is said that $\pmb{x}$ dominates $\pmb{y}$ (denoted as $\pmb{x} \succ \pmb{y}$) if $f_k(\pmb{x}) \leq f_k(\pmb{y})$ $ \forall k \in \{1,2\}$ and $f_k(\pmb{x})<f_k(\pmb{y})$ for at least one $k \in \{1,2\}$. When solving a biobjective problem we look for the Pareto front, which is the set composed by the non dominated solutions.

The exact Pareto frontier can be obtained using the improved version of the augmented $\epsilon$-constraint method (AUGMECON2), presented by \cite{Mavrotas2013}. 

According to this method, the problem needs to satisfy two conditions to determine the exact Pareto set: the objective function coefficients must be integer and the nadir points of the Pareto set must be known. The first condition can always be achieved by multiplying the coefficients of the objective functions by the appropriate power of 10. The second condition is met because we work with two objectives, which means that the nadir points can be obtained by calculating the pay-off table.
 
Basically, the AUGMECON2 method solves iteratively the following problem (denoted by $\hat{P}$):
{\allowdisplaybreaks
	\begin{alignat*}{2}
		\textrm{ min }  &
		[f_1(\pmb{x}) + \varepsilon (S_2/r_2)] \\ 
		\textrm{ st. } & f_2(\pmb{x}) - S_2 = e_2\\
		&  \pmb{x} \in F,
\end{alignat*}}

\noindent where $F$ represents the Constraints (\ref{constraint05}) - (\ref{constraint39}) of our problem, $S_2$ is the surplus variable for $f_2$, $ub_2$ is the upper bound of $f_2$, $r_2$ is the range of $f_2$, $g_2$ indicates the number of intervals that divide the range, $i_2$ is the grid point counter, $e_2 = ub_2 - (i_2 \times (r_2/g_2)) $ is the right hand side of the new constraint and $\varepsilon$ is a small number (usually between $10^{-3}$ and $10^{-6}$).

\begin{figure}[H]
	\begin{subfigure}{0.5\textwidth}
		\centering
		\includegraphics[scale=0.85,trim={0cm 0cm 0cm 0cm}, clip]{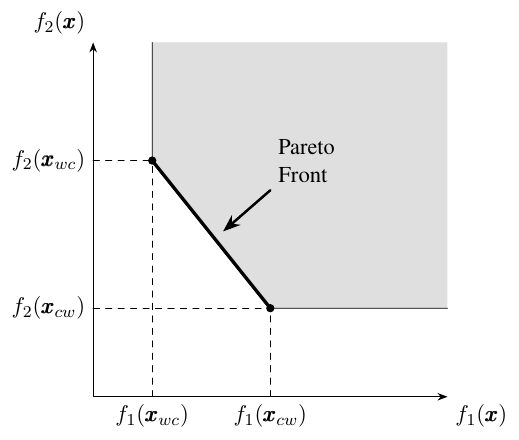}
		\label{fig:epsilon1}
	\end{subfigure}
	\begin{subfigure}{0.5\textwidth}
		\centering
		\includegraphics[scale=0.85,trim={0cm 0cm 0cm 0cm}, clip]{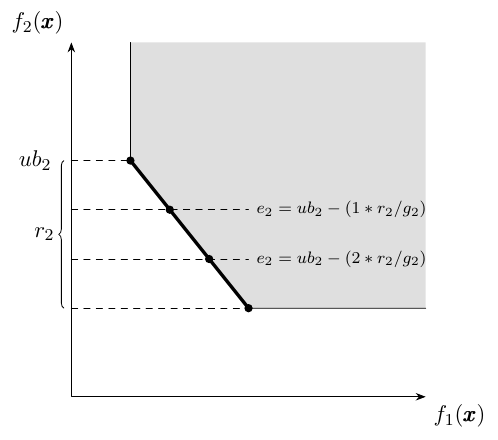}
		\label{fig:epsilon2}
	\end{subfigure}
	\caption{Example illustrating the AUGMECON2 method.}
	\label{fig:epsilon}
\end{figure}

To obtain $ub_2$ and $r_2$ it is necessary to solve the two lexicographical versions of our problem. It results in two solutions: one that prioritizes welfare over cost ($\pmb{x}_{wc}$) and another that prioritizes cost over welfare ($\pmb{x}_{cw}$). As it can be seen in Figure~\ref{fig:epsilon}, the upper bound is $ub_2=f_2(\pmb{x}_{wc})$ and the range is $r_2=f_2(\pmb{x}_{wc})-f_2(\pmb{x}_{cw})$. In the example of Figure~\ref{fig:epsilon} it is illustrated how the values of $e_2$ are obtained after dividing the range in $g_2 = 3$ intervals.

\vspace*{0.5cm}
\begin{algorithm}[H]
	\setstretch{1.05}
	\SetKwComment{Comment}{// }{}
	\caption{AUGMECON2 method} \label{alg:AUGMECON2}
	\KwData{Epsilon parameter ($\varepsilon$), range ($r_2$), upper bound ($ub_2$), number of intervals ($g_2$)}

	$i_2$ $\gets$ $1$ \\ 
	\While{$i_{2} \leq g_{2}$}{
		$\hat{P}$ $\gets$ \textit{{\textbf{generateMILP}}($ub_2, r_2, g_2, i_2, \varepsilon$)} \Comment{Generate the problem}
		$x$ $\gets$ \textit{{\textbf{solveMILP}}($\hat{P}$)} \Comment{Obtain a solution to the problem}
		\uIf{$x$ is feasible}{
			$\Omega$ $\gets$ $\Omega \cup \{\pmb{x}\}$ \Comment{Add the solution to the non dominated set}
			$b$ $\gets$ $\lfloor S_2/(r_2/g_2) \rfloor$, $i_2$ $\gets$ $i_2 + b$ \Comment{Get the number of iterations we can skip}
		}
		\uElse{
			\textbf{break} \Comment{The solution is not feasible}
		}
		\uIf{$i_2 < g_2$}{
			$i_2$ $\gets$ $i_2 + 1$ \Comment{Next grid point}
		}
	}
	
	\Return $\Omega$
	
\end{algorithm}
\vspace*{0.5cm}

The AUGMECON2 method, presented in Algorithm~\ref{alg:AUGMECON2}, initializes the grid point counters, $i_2$, (line 1) and iterates through them to solve the problem $\hat{P}$ (lines 3 - 4). If the solution found is feasible then it is added to the non dominated set (lines 5 - 6). The bypass coefficient, $b=	\lfloor S_2/(r_2/g_2)	\rfloor$, indicates the number of consecutive iterations that we can skip, $i_2 = i_2 + b$ (line 7). If the solution of $\hat{P}$ is not feasible (lines 8 - 9) or if we iterated through all the grid points (lines 10 - 11) the algorithm ends.

\subsection{The NSGA-II based algorithm}
In this part a Non-Dominated Sorting Genetic Algorithm II (NSGA-II)~\cite{Deb2002} is introduced, adapting it to our specific problem. To this aim, we divide the problem in two phases: on one hand, determining the routes of caregivers per day and, on the other hand, obtaining the optimal schedules for those routes. Note that the main decision variables involved in the routing phase are ($\var{x}{jk}$), while in the scheduling phase are the ones giving the starting times of the services ($\var{t}{j}$).

The solution encoding is really relevant when implementing this type of algorithms. For each day $d\in D$ and caregiver $i \in N$, we use a string vector to encode the routes assigned to each caregiver represented by $R^{id}=\{s_{i_{1}},s_{i_{2}},\ldots,s_{i_{d}}\}$, where $s_{i_{1}},s_{i_{2}},\ldots,s_{i_{d}}$ are the services IDs to be visited by the caregiver $i$ in the order from $s_{i_{1}}$ to $s_{i_{d}}$. Therefore, the routes of all caregivers in a day $d$ can be encoded by concatenating (using a unique identifier separator between the router of each caregiver) the string vectors of all caregivers in that day:
$$
R^{d} = \{R^{1d}+d_{1}+R^{2d}+d_{2}+\ldots+d_{n-1}+R^{nd}\}
$$
where $d_{1},d_{2},\ldots,d_{n-1}$ are unique separators to identify the routes associated to each caregiver.

In this way, the codification of the solution inside the NSGA-II is a chromosome divided into two parts represented by a tuple \texttt{<}$R,T$\texttt{>}, where $R$ is the string representing the routes of the solution previously explained, and $T$ contains the information about the schedule of each route. The general strategy we follow is:
\begin{itemize}
    \item The traditional operators of NSGA-II (crossover and mutation operators) are applied to the first part of the solution encoding related to the routes of the solutions, $R$.
    \item Once the routes are known, it is necessary to obtain a schedule $T$ associated to them in order to evaluate the objective function. Thus, the algorithms presented in Subsections~\ref{sec:imprwelf} and \ref{sec:imprcost} to find a schedule improving the welfare and the cost, respectively, are employed. It is important to note that giving a route, these algorithms can achieve several schedules. For this reason, it has been necessary to adapt NSGA-II so that, during the evaluation phase, all solutions obtained when calling those algorithms are incorporated into the existing population. By proceeding in this manner, it is possible to get solutions that aim to improve both objective functions at each iteration of the algorithm.
\end{itemize}
The procedure used to generate the initial solution is the same as the one employed for BIALNS. Specifically, the random greedy insertion operator is applied, either prioritizing welfare over cost or cost over welfare. Further, a uniform distribution in the selection of the initial population is applied to ensure that a similar number of individuals is created according to each one of the objectives of our problem.

\section{Computational experiments}

In this section we present the computational study carried out in order to check the behaviour of the biobjective algorithm described in the previous section.

Since there are two different objectives involved in the problem, it is not immediate to evaluate the quality of a solution or to decide which solutions are the best ones. Therefore, to compare the solutions of the biobjective problem, a set of performance indicators are used. The formal definition of the indicators is presented in Table~\ref{table:indicators} and they measure the convergence and the diversity of the solutions. To introduce them, it is necessary to consider an approximation of the Pareto frontier, $A$, and a reference set, $RF$. The set $A$ is the solution of the problem that has to be evaluated (for example, an approximation of the Pareto frontier obtained using BIALNS). Meanwhile, the reference set $RF$ is composed exclusively of non dominated points (since this set is usually not known a priori, a common approach is to define $RF$ by selecting the non dominated points of the solutions that are being evaluated by the indicator).

The indicators studied are the following ones. 

\begin{description}
	\item[Coverage.] It represents the percentage of elements of $A$ dominated by $RF$. The smaller the value of this measure, the better the quality of the approximation of the Pareto front.
	
	\item[Generational Distance.] It measures how far are the elements of $A$ from those of $RF$. This indicator is obtained using $d_i$, which is the euclidean distance between the elements $i \in A$ and the nearest one of $RF$. The smaller the value of this measure, the closer $A$ is to $RF$.
	
	\item[Inverted Generational Distance.] This variant of GD measures how far are the elements of the reference set $RF$ to the set $A$. Now, it takes the euclidean distance between the element $i \in RF$ and the nearest one of $A$, represented by $\tilde{d}_i$.  For this indicator, smaller values are preferred because it means that $RF$ is close to the approximation $A$.
	
	\item[Epsilon.]  It computes the minimum distance needed to translate every element of $A$ so it dominates the solution $RF$. It is said that $x 	\succ_{\epsilon}  y$ if, for each objective $k \in \{1,...,p\}$, $f_k(x) < \epsilon + f_k(y)$. If the value of this indicator is small, all the elements of solution $A$ are close to solution $RF$, because it is necessary to translate them a small distance in order to achieve the dominance of $RF$. The smaller the value, the better the quality of the approximation of the Pareto front.
\end{description}

\begin{table}[H]
	\renewcommand{\arraystretch}{1.1}
	\centering
	\begin{tabular}{l|c} 
		\textbf{Indicator}    &  \textbf{Formula}    \\
		\hline
		
		\rule{0pt}{21pt}Coverage (CV) &    \(\displaystyle  CV(RF,A)=   \frac{| \{ x \in A: \exists y \in RF / y \succ x\}|}{|A|}\)            \\[2ex] \hline
		
		\rule{0pt}{27pt}Generational Distance (GD) & \(\displaystyle  GD(RF,A)= \frac{\sqrt{\sum_{i=1}^{|A|} d_{i}^2}}{|A|}\) \\[2ex] \hline
		
		\rule{0pt}{27pt}Inverted Generational Distance (IGD) & \(\displaystyle   IGD(RF,A)= \frac{\sqrt{\sum_{i=1}^{|RF|} \tilde{d}_i^2}}{|RF|}  \)  \\[2ex] \hline
		
		\rule{0pt}{15pt}Epsilon (EPS) & \(\displaystyle    EPS(RF,A)=     \inf_{\epsilon \in \mathbb{R}} \{ \forall y \in RF \;\exists x \in A : x 	\succ_{\epsilon}  y  \}  \)  \\[1ex] \hline
	\end{tabular}
	\caption{Performance indicators for the biobjective problem.}
	\label{table:indicators}
\end{table}

All these indicators, except coverage, are implemented in the jMetal framework, described in \citet{DURILLO2011760}. The resolution approaches, BIALNS, AUGMECON2 and NSGA-II, were implemented in Python 3.7 (\citet{VanRossum2009}). Specifically, the NSGA-II algorithm has been implemented using the library jMetalPy~\cite{jmetalpy2019}. Furthermore, the MILP problem has been solved with Gurobi 9.1.1 (\citet{GurobiOptimization2021}) via its Python interface. All experiments were run in a machine Intel(R) Xeon(R) Gold 6146 CPU 3.20GHz, with 16GB of RAM, 2 cores and 100GB of hard drive, located in Centre for Information and Communications Technology Research (CITIC).

\subsection{Parameter study}

The analysis presented below is used to adjust the parameters involved in the algorithm. This is done using the instances  presented in \citet{Solomon1987}, which capture diverse scenarios in terms of service location and caregivers availability.  The instances specify duration, location and hard time windows of  services, as well as caregivers availability and maximum working time. Thus, the soft time windows of the services and the affinity levels between caregivers and services were randomly generated. To carry out the computational study, two types of Solomon instances were considered (with 10 and 15 services) and, for each of them, 10 different instances were randomly generated. 

\subsubsection{AUGMECON2 solutions}

The AUGMECON2 method was used to solve the Solomon instances with 10 and 15 services. For this purpose, a time limit of 12 hours has been established to solve each lexicographical MILP with the optimization solver Gurobi, in order to obtain the range of the objective functions necessary to define the grid points. Then, the MILP associated to each grid point has been solved with Gurobi with a time limit of 1 hour.  Table~\ref{table:AUGMECON2_time} presents the computational times, in hours, needed to solve each instance\footnote{Notice that, because of the time limit considered, there is no guarantee that the Pareto frontier of AUGMECON2 method is the optimal one.}. The fastest one is solved in 2.58 hours, while the slowest one needs 236.69 hours. Therefore, the AUGMECON2 method will not be suitable to solve large size instances. However, the results obtained with this method might be useful to evaluate the performance of the BIALNS algorithm.

\begin{table}[H]
	\centering
	\begin{tabular}{c|S[table-format=2.2]S[table-format=2.2]S[table-format=2.2]S[table-format=2.2]S[table-format=2.2]S[table-format=3.2]S[table-format=3.2]S[table-format=3.2]S[table-format=3.2]S[table-format=2.2]}
		Instance & \multicolumn{1}{c}{01} & \multicolumn{1}{c}{02} & \multicolumn{1}{c}{03} & \multicolumn{1}{c}{04} & \multicolumn{1}{c}{05} & \multicolumn{1}{c}{06} & \multicolumn{1}{c}{07} & \multicolumn{1}{c}{08} & \multicolumn{1}{c}{09} & \multicolumn{1}{c}{10} \\ \hline
		10       & 4.92                   & 2.58                   & 7.23                   & 20.39                  & 5.22                   & 5.83                   & 7.08                   & 5.13                   & 14.99                  & 3.58                    \\
		15       & 77.68                  & 36.67                  & 13.05                  & 58.97                  & 86.51                  & 236.69                 & 143.44                 & 124.00                 & 164.72                 & 90.29                  
		
	\end{tabular}  \caption{AUGMECON2 computational times (in hours).}\label{table:AUGMECON2_time}
\end{table}

\subsubsection{BIALNS parameter analysis}

This section focuses on studying BIALNS. The first step of the study is to evaluate the different parameters of the algorithm, in order to select the ones that find best solutions.

The algorithm is divided into three parts and, for each of them, different parameters were considered:

\begin{description}
	\setlength\itemsep{-0.25em}
	\item[Step 1: Initialize the sets.] In the first step it is necessary to establish the parameters for each lexicographic ALNS method. The values used are the ones obtained in \cite{MendezFernandez_thesis}:
	\begin{description}[topsep=-5.70pt]
		\setlength\itemsep{-0.25em}
		\item[ALNS.] Number of iterations (n): 1000. Proportion of solution to destroy (p): $auto\_100\%$.
	\end{description}
	\item[Step 2: Generate different solutions.] This step deals with the parameters related to the stopping criteria, as well as the parameters for each lexicographic ALNS method. 
	\begin{description}[topsep=-5.70pt]
		\setlength\itemsep{-0.25em}
		\item[Stopping criteria.] Number of iterations (nroutes): 1000, 2000, 3000, 4000, \textbf{5000}, \textbf{6000}, 7000, 8000, 9000 and 10000.
		\item[ALNS.] Number of iterations (nalns): \textbf{5}, \textbf{10}, 25, 50.  Proportion of solution to destroy (pr): \textbf{auto\_5\%}, \textbf{auto\_10\%} and auto\_25\%.
	\end{description}
	\item[Step 3: Get non dominated solutions.] In the third step, it is necessary to fix the stopping criteria parameter. 
	\begin{description}[topsep=-5.70pt]
		\item[Stopping criteria.] Number of iterations (nsols): $1\times10^{5}$, $\mathbf{2\times10^{5}}$ and $\mathbf{3\times10^{5}}$.
	\end{description}
\end{description}

The values highlighted in bold are the ones that provided best results (when comparing the solutions found by the algorithm with the ones obtained with AUGMECON2 method using performance indicators) for the instances with 10 and 15 services.
 
\begin{figure}[H]
	\centering
	\includegraphics[width=1\textwidth,trim={0cm 0cm 0cm 0cm}, clip]{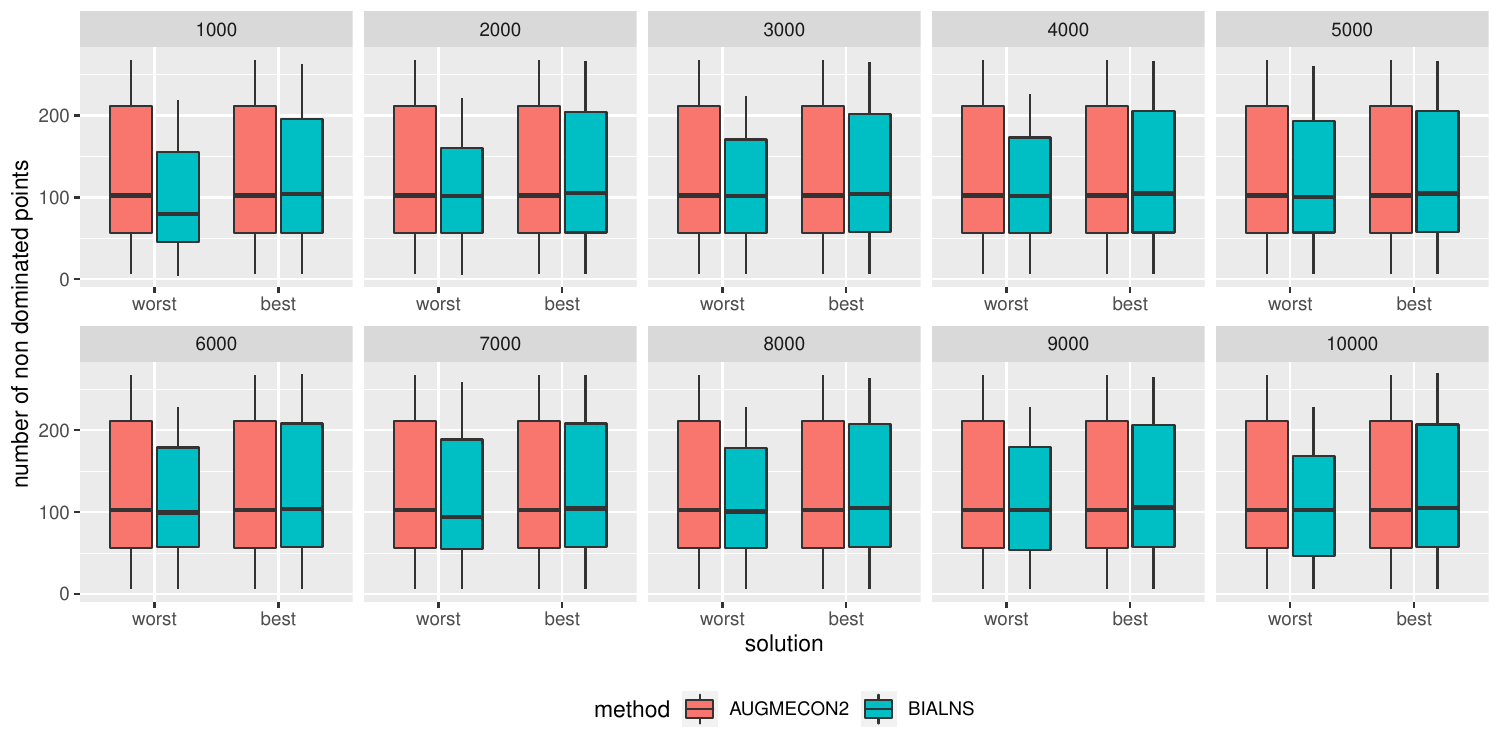}
	\caption{Number of non dominated points per iterations (10 services).}
	\label{fig:nd_points_10}
\end{figure}

\begin{figure}[H]
	\centering
	\includegraphics[width=1\textwidth,trim={0cm 0cm 0cm 0cm}, clip]{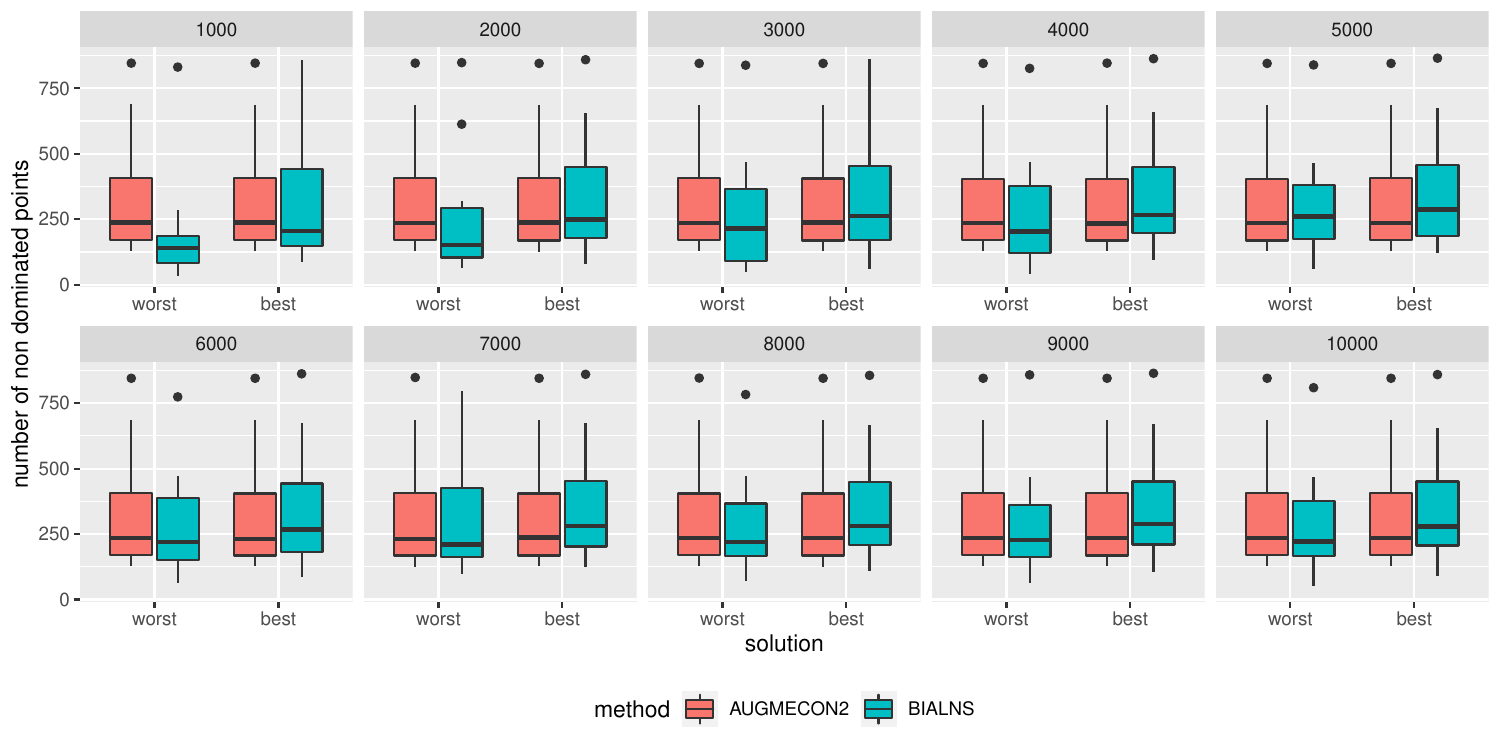}
	\caption{Number of non dominated points per iterations (15 services).}
	\label{fig:nd_points_15}
\end{figure}

To illustrate it, Figures~\ref{fig:nd_points_10}~and~\ref{fig:nd_points_15} show the number of non dominated points for the different values of the parameter nroutes. For  instances with 10 services it can be observed that, when the best solution is chosen, the algorithm is very competitive with AUGMECON2 in terms of the number of non dominated points. However, this is not the case of the worse solution. Specifically, the number of non dominated points for the worst solution increases with the iterations, reaching its best value from 2000 onward. In terms of the instances with 15 services, the best solution is practically not affected by the parameter, usually finding more non dominated points than AUGMECON2 method. Meanwhile, when considering a small number of iterations, the number of points is notably smaller for the worst solution. However, from 5000 iterations onward the worst solutions are similar to AUGMECON2 results.

Tables~\ref{tab:res_bialns_10services}~and~\ref{tab:res_bialns_15services} present the indicator values for the instances with 10 and 15 services with 6000 and 8000 iterations, respectively. The number of iterations considered is chosen according to the results presented in Figures~\ref{fig:nd_points_10}~and~\ref{fig:nd_points_15}, but a more detailed analysis can be found in \citet{MendezFernandez_thesis}, where confidence intervals are used. In terms of the instances with 10 services, Table~\ref{tab:res_bialns_10services} shows that no method performs better than the other. Instead, the behaviour of the methods is instance dependant. Regrading the instances with 15 services, CV values indicate that the proportion of dominated points of the algorithm is larger than the one of the AUGMECON2 method. Meanwhile, indicators EPS, GD and IGD show that BIALNS solutions are better than the ones of AUGMECON2.

\begin{table}[htbp]
  \centering
  \begin{tabular}{|l|l|l|l|l|l|l|l|l|}
    \hline
    Instances & \multicolumn{4}{c|}{AUGMECON2} & \multicolumn{4}{c|}{BIALNS} \\ \hline
              & CV           & EPS           &  GD           &     IGD       & CV            & EPS           &  GD           &     IGD       \\ \hline
    10\_01 &  0 &                    0 &       0 &                    0 &  2.1e-2 &  3.8e-3 & 5.8e{-5} &              1.2e-4             \\ \hline
    10\_02 & 0 & 2.6e{-5} &       0 & 2.4e{-5} &  7.4e-2 &  2.3e-3 & 1.6e{-5} &              1.2e-4          \\ \hline
    10\_03 & 0 & 5.2e{-5} &       0 &              3.9e-4 &      0 &       0 &                    0 &                    0           \\ \hline
    10\_04 & 8.3e-2 &                6.5e-2 &   1.6e-2 &                3.1e-2 &      0 &       0 &                    0 &                    0       \\ \hline
    10\_05 & 0 & 4.7e{-5} &       0 &              1.6e-4 &      0 &       0 &                    0 &                    0 \\ \hline
    10\_06 &   0 &                    0 &       0 &                    0 &  4.4e-2 &   1.4e-2 &              5.1e-4 &              5.4e-4  \\ \hline
    10\_07 &    0 &              1.5e-4 &       0 &              1.4e-4 &  4.8e-2 &  8.4e-3 &              1.7e-4 &               3e-4       \\ \hline
    10\_08 &   0 &                    0 &       0 &                    0 &  4.6e-2 &  3.3e-3& 6.1e{-5}& 9.2e{-5}       \\ \hline
    10\_09 & 1.3e-2 &                1.3e-2 & 2.1e-4 &              6.3e-4 &  2.1e-2 &  7.6e-3 & 9.6e{-5} &              1.4e-4         \\ \hline
    10\_10 &   0 &                    0 &       0 &                    0 &      0 &       0 &                    0 &                    0      \\ \hline
  \end{tabular}
	\caption{Comparison of AUGMECON2 and BIALNS (10 services).}
  \label{tab:res_bialns_10services}
\end{table}

\begin{table}[htbp]
  \centering
  \begin{tabular}{|l|l|l|l|l|l|l|l|l|}
    \hline
    Instances & \multicolumn{4}{c|}{AUGMECON2} & \multicolumn{4}{c|}{BIALNS} \\ \hline
              & CV           & EPS           &  GD           &     IGD       & CV            & EPS           &  GD           &     IGD       \\ \hline
    15\_01 & 1.0e-2 &               1.4e-3 &             1.0e{-6} &             2.9e{-5} &   1.6e-1 &  4.4e-3 &              3.2e-4 &              1.1e-4  \\ \hline
    15\_02& 1.7e-2 &                2.8e-2&             9.7e-3 &               1.5e-3 & 4.6e-3 &  1.2e-3 &             4.7e{-6} &               1.6e-4       \\ \hline
    15\_03 &    0 &             3.8e{-5} &                  0 &             9.1e{-5} &      0 & 1.5e-4&                    0 &              1.2e-4    \\ \hline
    15\_04 & 5.4e-3 &               2.2e-3 &           3.2e-4&              3.5e-4 &  1.8e-2 &   3.0e-3 &             1.5e{-5} & 6.5e{-5}        \\ \hline
    15\_05 & 3.9e-2 &                4.3e-2 &             2.9e-3&               2.8e-3 &   1.6e-1 &  3.4e-3 &              3.9e-4&              2.8e-4        \\ \hline
    15\_06 &7.1e-3 &                6.8e-2 &             1.9e-3 &                1.1e-2& 4.7e-3 &  2.2e-3 &             2.2e{-5} & 1.1e{-5}\\ \hline
    15\_07 &   6.5e-2 &                 1.2e-1&             1.9e-3 &                 2.0e-2&   2.6e-1 &  2.5e-3 &             6.3e-4 &               6.0e-4   \\ \hline
    15\_08 & 3.6e-2&               5.3e-2 &            7.7e-4 &                2.6e-2 &   2.5e-1 &  4.1e-3 &              3.2e-4 &               1.4e-3   \\ \hline
    15\_09 &   2.0e-2&                 1.3e-1&             3.9e-3&                5.9e-2 &   5.8e-1 &  4.3e-3&               1.9e-3&               2.7e-3    \\ \hline
    15\_10 &  2.5e-2 &               9.7e-3&             1.5e-3 &              2.6e-4 &   2.6e-1&  1.1e-3 &             2.4e{-5} & 6.8e{-5}         \\ \hline
  \end{tabular}
	\caption{Comparison of AUGMECON2 and BIALNS (15 services).}
  \label{tab:res_bialns_15services}
\end{table}

\begin{figure}[H]
	\begin{subfigure}{0.5\textwidth}
		\centering
		\includegraphics[width=0.99\textwidth,trim={0cm 0cm 0cm 0cm}, clip]{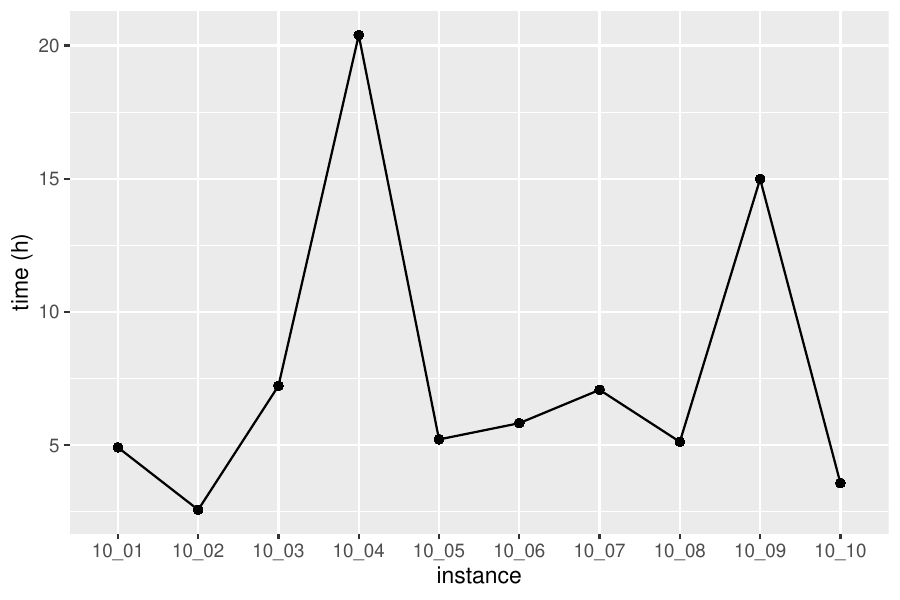}
		\caption{AUGMECON2 time (hours).}
		\label{fig:eps_time_10}
	\end{subfigure}
	\begin{subfigure}{0.5\textwidth}
		\centering
		\includegraphics[width=0.99\textwidth,trim={0cm 0cm 0cm 0cm}, clip]{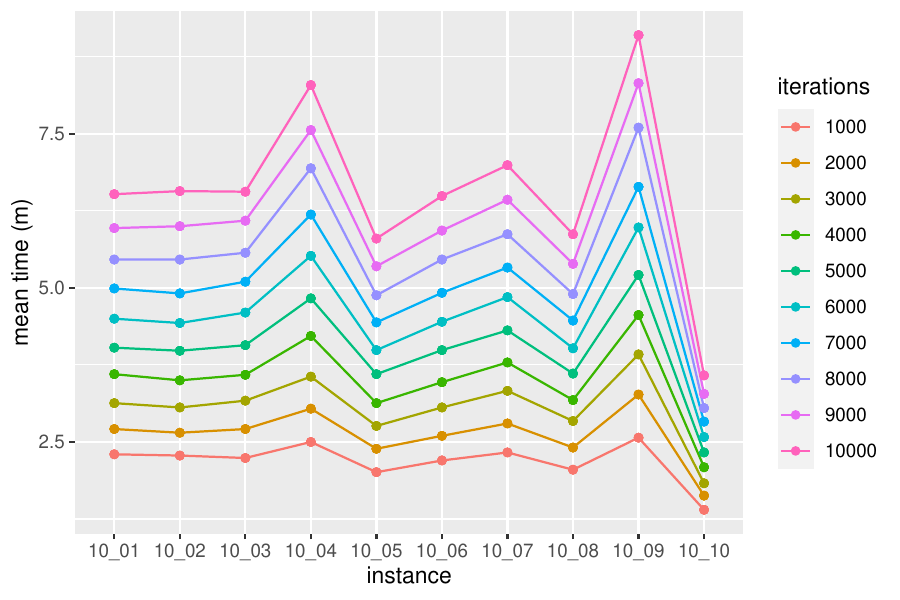}
		\caption{BIALNS time (minutes).}
		\label{fig:alg_time_10}
	\end{subfigure}
	\caption{Computational times (10 services).}
	\label{fig:times_10}
\end{figure}

\begin{figure}[H]
	\begin{subfigure}{0.5\textwidth}
		\centering
		\includegraphics[width=0.99\textwidth,trim={0cm 0cm 0cm 0cm}, clip]{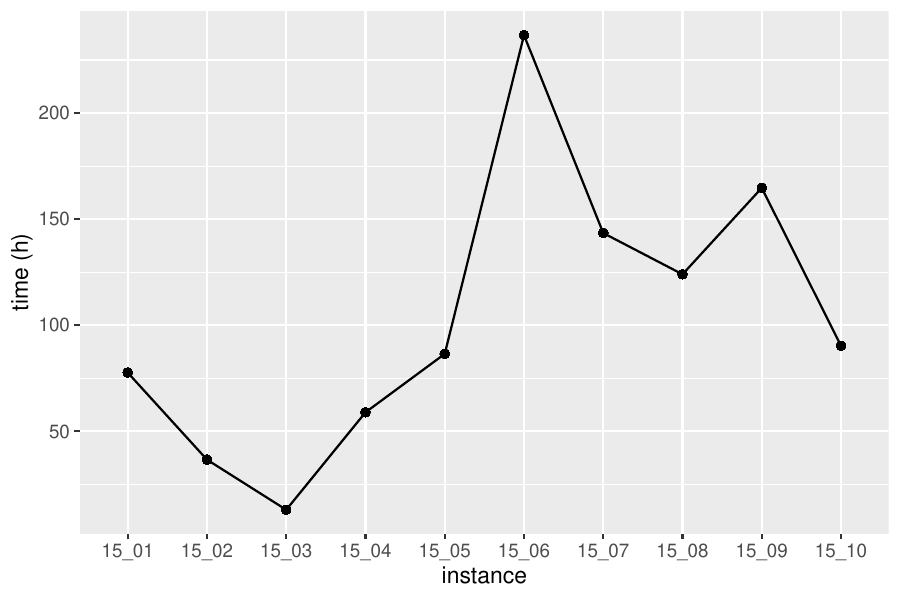}
		\caption{AUGMECON2 time (hours).}
		\label{fig:eps_time_15}
	\end{subfigure}
	\begin{subfigure}{0.5\textwidth}
		\centering
		\includegraphics[width=0.99\textwidth,trim={0cm 0cm 0cm 0cm}, clip]{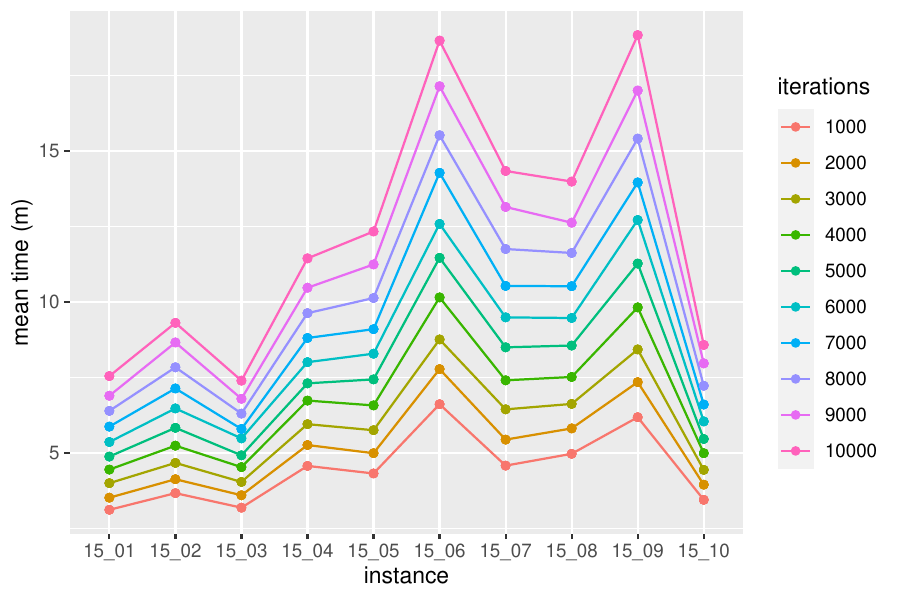}
		\caption{BIALNS time (minutes).}
		\label{fig:alg_time_15}
	\end{subfigure}
	\caption{Computational times (15 services).}
	\label{fig:times_15}
\end{figure}

Figures~\ref{fig:times_10}~and~\ref{fig:times_15} show computational times\footnote{Notice that the algorithm was run 5 times. Therefore, the figure presents the mean time employed to solve each instance for each value of the parameter.} after solving the instances with 10 and 15 services, respectively.  Notice that, in both methods, the behavior of the instances in terms of computational time is similar but in different magnitudes.  For the 10 service instances, the AUGMECON2 method needs more than 4 hours to solve most of the considered instances, meanwhile the BIALNS solves all of them in less than 10 minutes. Interestingly, the comparison of the two figures shows that the most time-consuming instances for the AUGMECON2 method ($10\_04$ and $10\_09$) are also the slowest ones for the algorithm. Similarly, in case of the 15 service instances, the computational times for the AUGMECON2 vary from 13.05 to 236.69 hours, meanwhile the highest computational time is lower than 20 minutes.

\subsubsection{Comparison with \sc{NSGA-II}}
The parameters of the NSGA-II algorithm have been selected after some preliminary experiments. In particular, the configuration employed is:
\begin{itemize}
	\item Population size: 250.
	\item Generations size: 500.
	\item Crossover operator: Partially Mapped Crossover (PMX) with probability 0.8.
	\item Mutation operator: Permutation Swap Mutation with probability 0.3.
	\item Selection operator: Binary Tournament Selection.
\end{itemize}

Table~\ref{tab:res_nsga_10services} reports  performance indicators comparing NSGA-II with BIALNS\footnote{BIALNS is run employing the best configuration parameters obtained after the analysis showed in previous section.} over  instances of 10 services. As in the case of BIALNS, NSGA-II was run five times per instance, reporting mean values for the metrics. Furthermore, a time limit of 600 seconds is considered as a stopping criterion for NSGA-II (a time limit greater than the slowest instance solved by BIALNS). Table~\ref{tab:res_nsga_10services} clearly shows that BIALNS outperforms NSGA-II according to the four metrics presented for all the instances. Considering CV, it shows that practically all  solutions found by NSGA-II are dominated by the ones found by BIALNS. Analysing remaining metrics, it can be seen that the values for NSGA-II are larger than those for BIALNS (it practically are always 0). It is worth highlighting instance R202\_025, in which NSGA-II exhibits its best performance, achieving the lowest metric values among the instances, although still larger than those of BIALNS. 
\begin{table}[htbp]
  \centering
  \begin{tabular}{|l|l|l|l|l|l|l|l|l|}
    \hline
    Instances & \multicolumn{4}{c|}{NSGA-II} & \multicolumn{4}{c|}{BIALNS} \\ \hline
              & CV           & EPS           &  GD           &     IGD       & CV            & EPS           &  GD           &     IGD       \\ \hline
    10\_01 & 0.8310       & 0.1705        & 0.0920        & 0.1325        & 0             & 0             & 0             & 0             \\ \hline
    10\_02 & 0.9800       & 0.0749        & 0.2516        & 0.1093        & 0             & 0             & 0             & 0             \\ \hline
    10\_03 & 1.0000       & 0.0465        & 0.2949        & 0.0864        & 0             & 0             & 0             & 0             \\ \hline
    10\_04 & 1.0000       & 0.3567        & 0.2160        & 0.2524        & 0             & 0             & 0             & 0             \\ \hline
    10\_05 & 0.8159       & 0.1508        & 0.0832        & 0.0643        & 0.0784        & 0.0082        & 0.0043        & 0.0063        \\ \hline
    10\_06 & 0.9778       & 0.2123        & 0.1252        & 0.1205        & 0             & 0             & 0             & 0             \\ \hline
    10\_07 & 0.9314       & 0.1750        & 0.2380        & 0.1579        & 0             & 0             & 0             & 0             \\ \hline
    10\_08 & 0.8776       & 0.0947        & 0.0743        & 0.0582        & 0             & 0             & 0             & 0             \\ \hline
    10\_09 & 1.0000       & 0.2204        & 0.2144        & 0.2057        & 0             & 0             & 0             & 0             \\ \hline
    10\_10 & 1.0000      & 0.6083        & 0.2310        & 0.4824        & 0             & 0             & 0             & 0             \\ \hline
  \end{tabular}
	\caption{Comparison of NSGA-II and BIALNS (10 services).}
  \label{tab:res_nsga_10services}
\end{table}

\begin{figure}[htbp]
    \begin{center}
		\includegraphics[trim=0 0 175 90,clip,width=0.95\textwidth]{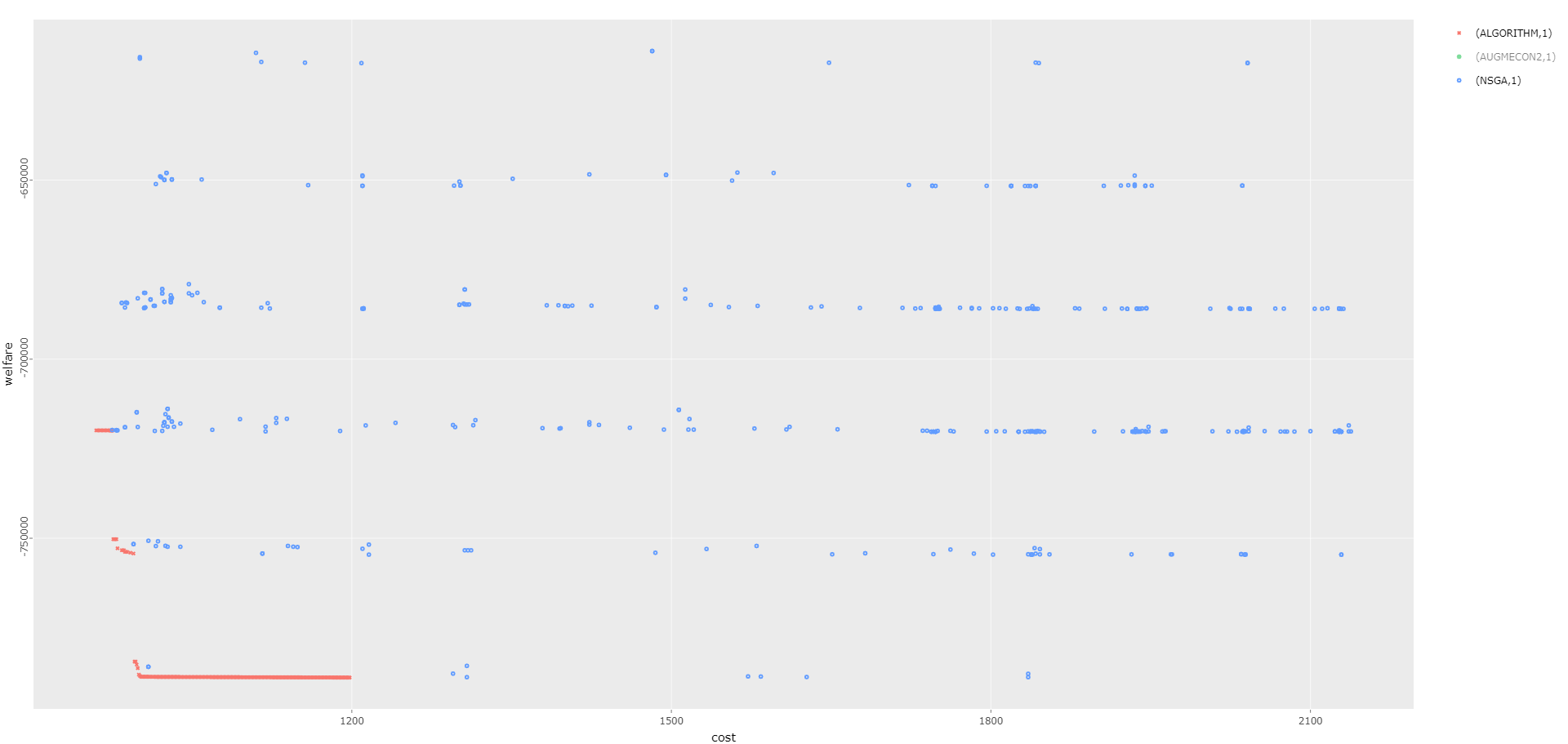}
    \end{center}
		\caption{Approximations of the Pareto fronts of BIALNS (red points) and NSGA-II (blue points) (10 services).}
    \label{fig:10_nsga_paretofront}
\end{figure}

In order to get a deeper understanding, Figure~\ref{fig:10_nsga_paretofront} represents the approximation of the Pareto front obtained with BIALNS and NSGA-II for one of the instances of 10 services\footnote{The behaviour for the remaining instances is similar.}. It is clear that BIALNS obtains a better set of points than NSGA-II. Although NSGA-II finds more scattered points (it captures more ``jumps'' corresponding to changes in the welfare, due to variations in the affinity levels caregiver-services), practically all of them are dominated by the points obtained by NSGA-II. These observations are in accordance with the performance metrics results discussed above.

The results for the instances with 15 services are similar, so for the sake of clarity, they are not included. 

\subsection{Real instances}

The real instances were obtained during consecutive weeks between 2016 and 2017. Next, the behaviour of BIALNS will be shown in one of these instances. The instance deals with 808 services to be carried out throughout a week by the 36 available caregivers. To solve the instance, next values of the parameters were taken into account:

\begin{description}
	\item[Step 1: Initialize the sets.] Proportion of solution to destroy:  $auto\_1\%$. Time limit: 90 minutes.
	\item[Step 2: Generate different solutions.]  Number of iterations (nroutes): 10000. Number of iterations (nalns): 1. Proportion of solution to destroy (pr): $1\%$.
	\item[Step 3: Get non dominated solutions.]  Number of iterations (nsols): $3\times10^{5}$.
\end{description}

\begin{figure}[H]
	\centering
	\includegraphics[width=1\textwidth,trim={0cm 0cm 0cm 0cm}, clip]{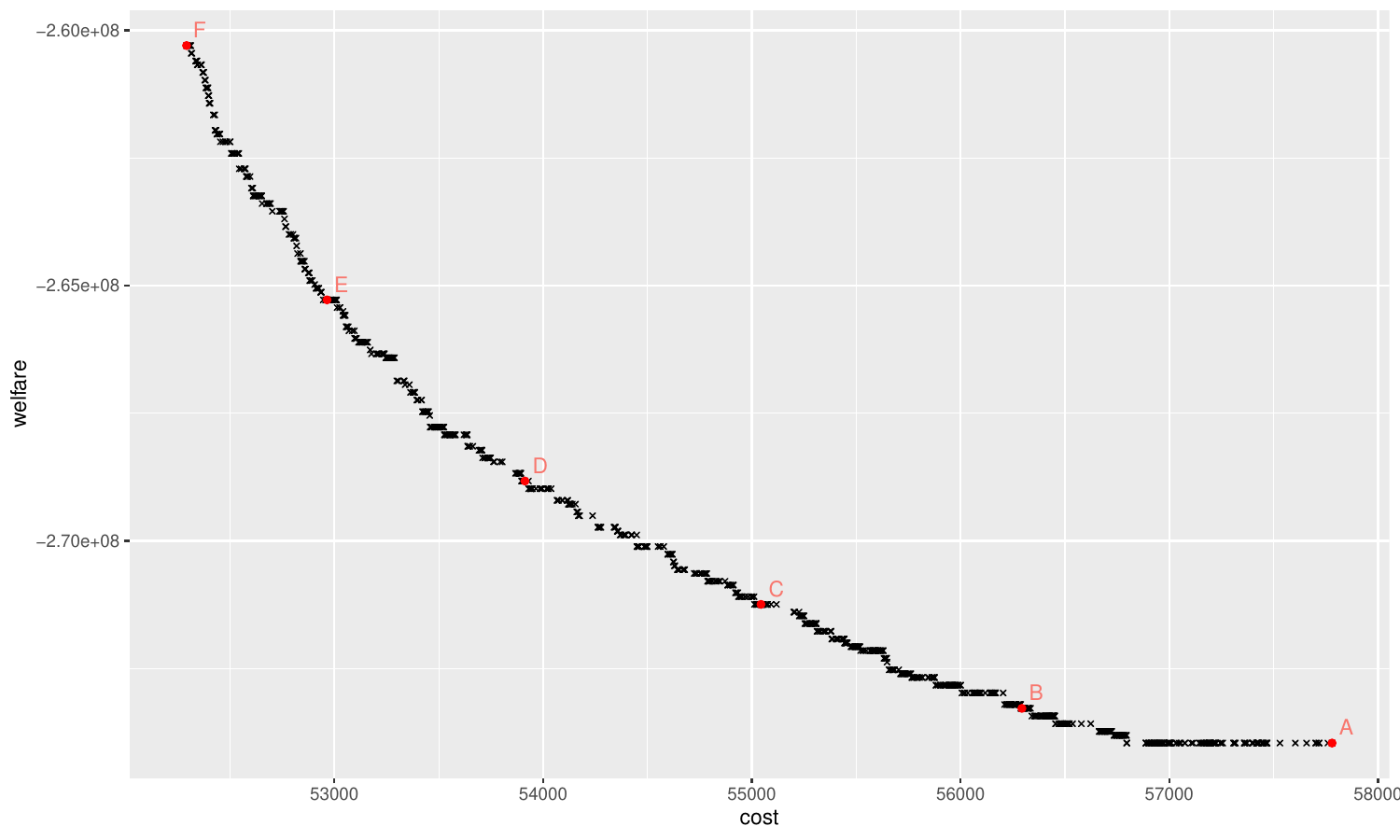}
	\caption{Pareto frontier for  week 9.}
	\label{fig:front_week9}
\end{figure}

Figure~\ref{fig:front_week9} shows the Pareto frontier for week 9. To analyze the solutions found, six points along the frontier have been selected (denoted as A - F). Point A (F) is the solution of the hierarchical problem that prioritizes welfare over cost ( cost over  welfare). Points B, C, D and E are intermediate solutions. 

Since all points presented in the frontier are non dominated, to improve welfare it is necessary to increase cost. This can be seen in more detail in Table~\ref{table:objetives_biobj}, which presents the objective values of the selected points, as well as the ones of the solution considered by the company. Regarding cost, all the solutions considered are better than the one of the company, which can also be assessed  analyzing the overtime and working time of the solutions. In terms of  welfare, the solution of the company is between points B and C, although the penalization value of C is much better than the one of the company.

\begin{table}[H]
	\centering
	\begin{tabular}{ccccccc}
		Solution & Cost  & Overtime & Worked time & Welfare ($\times e{8}$) & Affinity & Penalization \\ \hline
		Company  & 69744 &  11525   &    58219    &          -2.71          &   3576   &    10518     \\
		   A     & 57780 &   4692   &    53088    &          -2.74          &   3631   &     3342     \\
		   B     & 56294 &   3761   &    52533    &          -2.73          &   3622   &     4005     \\
		   C     & 55044 &   2733   &    52311    &          -2.71          &   3595   &     4436     \\
		   D     & 53913 &   1670   &    52243    &          -2.68          &   3563   &     5218     \\
		   E     & 52965 &   722    &    52243    &          -2.65          &   3516   &     7900     \\
		   F     & 52292 &   273    &    52019    &          -2.60          &   3450   &     9961     \\ \hline
	\end{tabular}
	\caption{Objective values of the solutions.}\label{table:objetives_biobj}
\end{table}

\begin{figure}[H]
	\centering
	\includegraphics[width=0.6\textwidth,trim={0cm 0cm 0cm 0cm}, clip]{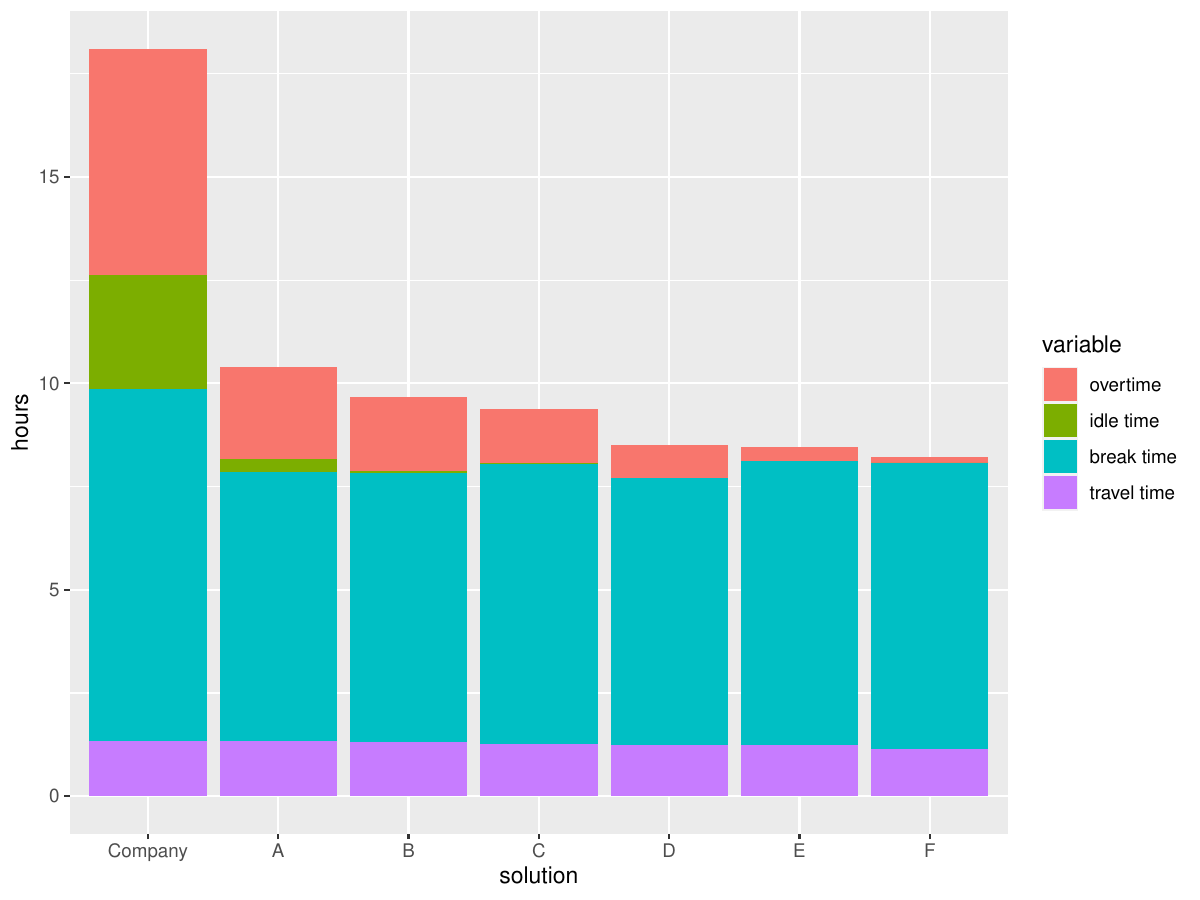}
	\caption{Mean weekly time per caregiver  (week 9).}
	\label{fig:cost_week9}
\end{figure}

\begin{figure}[H]
	\centering
	\includegraphics[width=0.6\textwidth,trim={0cm 0cm 0cm 0cm}, clip]{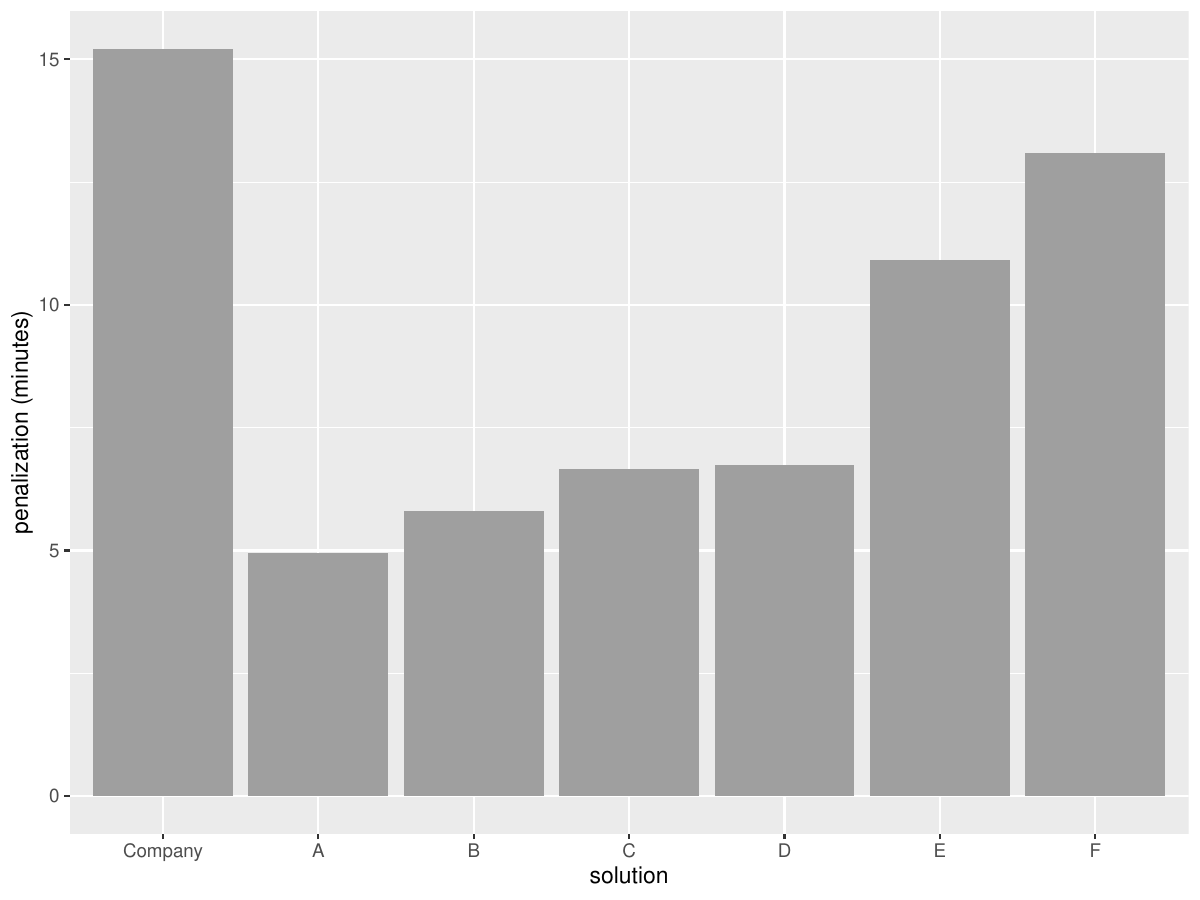}
	\caption{Mean soft time window penalization per user (week 9).}
	\label{fig:stw_week9}
\end{figure}

To complete the comparison,  Figures~\ref{fig:cost_week9}~and~\ref{fig:stw_week9} are presented.
Figure~\ref{fig:cost_week9} represents the mean overtime, idle time (i.e. paid break time), break time (i.e. unpaid break time) and travel time per caregiver according to the considered solutions. Figure~\ref{fig:stw_week9} displays the mean soft time window penalization per user. According to these figures it is clear that the six solutions taken from the Pareto frontier are better than the one of the company. Specifically, it can be seen that solution A (which prioritizes welfare over  cost) is the one with highest overtime and idle time but, as more priority is given to cost,  idle time disappears and overtime decreases. In a similar way, Figure~\ref{fig:stw_week9} shows that as the priority of welfare decreases, soft time window penalization increases.

\section{Concluding remarks}

In this work, a real problem of a home care company is tackled by considering it as a biobjective HCSP. The problem shares many of the characteristics of other routing and scheduling multiobjective problems in home care. The work was motivated by the company's interest in finding a solution that balances two objectives: the welfare of users and caregivers, and the cost associated with the planning. Furthermore, there is a special feature that substantially distinguishes it from the others and makes its study of great interest: the longest break between two daily consecutive services of each caregiver will not be included as part of her working day, provided that it is greater than a fixed number of hours.

The problem is approached from several perspectives. First, a metaheuristic MDLS algorithm, called BIALNS, is proposed. BIALNS involves generating different solutions and then modifying the schedules to obtain non-dominated solutions using an ALNS technique and two heuristic scheduling methods designed for hierarchical problems. BIALNS is compared with other two well-known multiobjective techniques, which were adapted to the problem specific requirements: AUGMECON2 and NSGA-II. The idea behind the implementation of AUGMECON2 was to obtain the exact Pareto front, to guarantee that the BIALNS algorithm provides good approximation in reasonable computational times. The widespread used of NSGA-II in multiobjective problems, makes it a good choice for comparison with BIALNS.

To validate the algorithm, several computational experiments using literature instances were run, considering different configurations for BIALNS algorithm, in order to select the ones that provide the best solutions. The quality of the solutions was analyzed by comparing them to the ones obtained with the AUGMECON2 and NSGA-II methods. After thorough experimentation, we can assert that our approach yields solutions of superior quality (based on common multiobjective performance metrics) compared to those generated by NSGA-II. On the other hand, AUGMECON2 managed to attain the exact Pareto frontier in some small instances, despite the high computational costs. In those cases, our approach was able to achieve very good approximations much more rapidly.

The algorithm presented in this work is devoted to help the company to make quick decisions about the routes and schedules of the caregivers of the company. For this reason, BIALNS is applied over a real instance of the company, showing several configurations for the schedules of the caregivers in terms of a trade-off between the welfare of users/caregivers and the cost of the schedules. Due to the good behaviour of BIALNS, it would be in the future a good choice to improve the expert system of the company.

\section*{CRediT authorship contribution statement}
\textbf{Isabel M\'{e}ndez-Fern\'{a}ndez:} Conceptualization, Methodology, Software, Writing – original draft, Writing – review \& editing. \textbf{Silvia Lorenzo-Freire:} Conceptualization, Methodology, Software, Writing – original draft, Writing – review \& editing. \textbf{\'{A}ngel Manuel Gonz\'{a}lez-Rueda:} Conceptualization, Methodology, Software, Writing – original draft, Writing – review \& editing.

\section*{Acknowledgements}
This research was funded by MICINN/AEI/10.13039/501100011033/ and ERDF/EU through R+D+I project grants  MTM2017-87197-C3-1-P and PID2021-124030NB-C31. Isabel M\'{e}ndez-Fern\'{a}ndez and Silvia Lorenzo-Freire also acknowledge support from Grupos de Referencia Competitiva ED431C-2020/14 and Centro de Investigaci\'{o}n del Sistema universitario de Galicia ED431G- 2019/01; and \'{A}ngel M. Gonz\'{a}lez-Rueda from Grupos de Referencia Competitiva ED431C-2021/24, all of them funded by Conseller\'{\i}a de Cultura, Educaci\'on e Universidades, Xunta de Galicia.

\bibliographystyle{elsarticle-num-names} 
\bibliography{biliography_mayores.bib}

\end{document}